\documentclass[11pt,leqno]{amsart}
\usepackage{amsmath, amssymb, amsfonts, amsthm}
\usepackage{xspace, amscd, graphics}

\theoremstyle{plain}

\newtheorem{main}{Theorem}

\newtheorem{theorem}{Theorem}[section]

\newtheorem{proposition}[theorem]{Proposition}
\newtheorem{corollary}[theorem]{Corollary}
\newtheorem{lemma}[theorem]{Lemma}

\newcommand{\Set}[1]{{\left\{\,#1\,\right\}}}
\newcommand{\Q}{{\mathbb Q}}
\newcommand{\Z}{{\mathbb Z}}
\newcommand{\R}{{\mathbb R}}
\newcommand{\C}{{\mathbb C}}
\newcommand{\G}{{\mathcal G}}

\newcommand{\A}{{\mathbb A}}

\begin{document}


\title{A Cuspidality Criterion for the Functorial Product on ${\bf GL
(2) \times GL (3)}$, with a cohomological application}
\author{Dinakar Ramakrishnan \and Song Wang}
\date{}
\maketitle

\pagestyle{myheadings} \markboth{Dinakar Ramakrishnan and Song
Wang}{A Cuspidality Criterion with a cohomological application}

\bigskip


\section{Introduction}

A strong impetus for this paper came, at least for the first
author, from a question of Avner Ash, asking whether one can
construct non-selfdual, non-monomial cuspidal cohomology classes
for suitable congruence subgroups $\Gamma$ of SL$(n, \Z)$, say for
$n=6$. Such a construction, in special examples, has been known
for some time for $n=3$ (\cite{A-G-G}, \cite{vG-T1},
\cite{vG-K-T-V}, \cite{vG-T2}); it is of course not possible for
$n=2$. One can without trouble construct non--selfdual,
\emph{monomial} classes for any $n = 2 m$ with $m \geq 2$, not
just for constant coefficients (see the Appendix, Theorem E). In
the Appendix we also construct non-monomial, non-selfdual classes
for $n=4$ using the automorphic induction to $\Q$ of suitable
Hecke character twists of non-CM cusp forms of ``weight $2$'' over
imaginary quadratic fields, but they admit quadratic self-twists
and are hence imprimitive. The tack pursued in the main body of
this paper, and which is the natural thing to do, is to take a
non-selfdual (non-monomial) $n=3$ example $\pi$, and take its
functorial product $\boxtimes$ with a cuspidal $\pi'$ on
GL$(2)/\Q$ associated to a holomorphic newform of weight $4$ for a
congruence subgroup of SL$(2, \Z)$. The resulting (cohomological)
$n=6$ example can be shown to be {\it non-selfdual} for suitable
$\pi'$. (This should be the case for all $\pi'$, but we cannot
prove this with current technology -- see Remark 4.1.) Given that,
the main problem is that it is not easy to show that such an
automorphic tensor product $\Pi:= \pi \boxtimes \pi'$, whose
modularity was established in the recent deep work of H.~Kim and
F.~Shahidi (\cite{KSh2000}), is {\it cuspidal}. This has led us to
prove a precise cuspidality criterion (Theorem A) for this
product, not just for those of cohomological type, which hopefully
justifies the existence of this paper. The second author earlier
proved such a criterion when $\pi$ is a twist of the symmetric
square of a cusp form on GL$(2)$ (\cite{Wa2003}; such forms are
essentially selfdual, however, and so do not help towards the
problem of constructing non-selfdual classes. One of the reasons
we are able to prove the criterion {\it in general} is the fact
that the associated, degree $20$ exterior cube $L$-function is
nicely behaved and analyzable. This helps us rule out, when the
forms on GL$(2)$ and GL$(3)$ are non-monomial, the possible
decomposition of $\Pi$ into an isobaric sum of two cusp forms on
GL$(3)$ (see section 7). This is the heart of the matter.

We will also give a criterion (Theorem B) as to when the base
change of $\Pi$ to a solvable Galois extension remains cuspidal.
We will derive a stronger result for the cohomological examples
(Theorem C), namely that each of them is {\it primitive}, i.e.,
not associated to a cusp form on GL$(m)/K$ for {\it any} (possibly
non-normal) extension $K/\Q$ of degree $d > 1$ with $dm=6$.
Furthermore, each of the three main non-selfdual GL$(3)$ examples
$\pi$ of \cite{vG-T1}, \cite{vG-K-T-V} and \cite{vG-T2} comes
equipped, confirming a basic conjecture of Clozel (\cite{C}), with
a certain $3$-dimensional $\ell$-adic representation $\rho_\ell$
whose Frobenius traces $a_p(\rho_\ell)$ agree with the Hecke
eigenvalues $a_p(\pi)$ for small $p$. For $\pi'$ on GL$(2)/\Q$
defined by a suitable holomorphic newform of weight $4$, with
associated Galois representation $\rho'_\ell$, we will show
(Theorem D) that the six-dimensional $R_\ell: = \rho_\ell \otimes
\rho'_\ell$, which should correspond to $\Pi$, remains irreducible
under restriction to {\it any} open subgroup of Gal$(\overline
\Q/\Q)$.

The first author would like to thank Avner Ash for his question
and for making comments on a earlier version of this paper, Mahdi
Asgari for initially kindling his interest (at Park City, UT) in
the problem of establishing a precise {\it cuspidality criterion}
for the Kim-Shahidi product, and the National Science Foundation
for financial support through the grant DMS-0100372. The second
author would like to thank James Cogdell and Henry Kim for their
interest in his lecture on this work at the Fields Institute
Workshop on Automorphic L-functions in May 2003.

\bigskip


\section{The Cuspidality Criterion}

Throughout this paper, by a cusp form on $GL (n)$ (over a global
field $F$) we will mean an irreducible, cuspidal automorphic
representation $\pi = \pi_\infty \otimes \pi_f$ of $GL_{n}
(\mathbb{A}_{F})$. We will denote its central character by
$\omega_\pi$. One says that $\pi$ is {\it essentially self-dual}
iff its contragredient $\pi^\vee$ is isomorphic to $\pi \otimes
\nu$ for some character $\nu$ of (the idele classes of) $F$; when
$n=2$, one always has $\pi^\vee \simeq \pi \otimes
\omega_\pi^{-1}$. Note that $\pi$ is {\it unitary} iff $\pi^\vee$
is the {\it complex conjugate representation} $\overline \pi$.
Given any cusp form $\pi$, we can find a real number $t$ such that
$\pi_u: = \pi \otimes \vert . \vert^t$ is unitary.

For any cusp form $\pi'$ on $GL (2)$, put $Ad (\pi') = {\rm
sym}^{2} (\pi') \otimes \omega_{\pi'}^{-1}$ and $A^{4} (\pi') =
{\rm sym}^{4} (\pi') \otimes \omega_{\pi'}^{-2}$. Recall that
$\pi'$ is {\it dihedral} iff it admits a self-twist by a quadratic
character; it is {\it tetrahedral}, resp. {\it octahedral}, iff
$sym^2(\pi')$, resp. $sym^3(\pi')$, is cuspidal and admits a
self-twist by a cubic, resp. quadratic, character. (The automorphy
of $sym^3(\pi')$ was shown by Kim and Shahidi in \cite{KSh2000}.)
We will say that $\pi'$ is of {\it solvable polyhedral type} iff
it is dihedral, tetrahedral or octahedral.

\begin{main} \label{TM:A}
    Let $\pi', \pi$ be cusp forms on $GL (2)$, $GL (3)$
    respectively over a number field $F$. Then the Kim--Shahidi transfer $\Pi =
        \pi \boxtimes \pi'$ on GL$(6)/F$
    is cuspidal unless one of the following happens:

        \textnormal{(a)} $\pi'$ is not dihedral, and $\pi$ is a twist of $Ad (\pi')$;

        \textnormal{(b)} $\pi'$ is dihedral, $L(s, \pi) = L(s,
        \chi)$ for an idele class character $\chi$ of a cubic,
        non-normal extension $K$ of $F$, and the base change $\pi_K$ is
        Eisensteinian.

    \medskip

    Furthermore, when \textnormal{(a)} (resp. \textnormal{(b)}) happens,
    $\Pi$ possesses an isobaric
    decomposition of type $(2, 4)$ or $(2, 2, 2)$  (resp. of type $(3, 3)$). More
    precisely, when we are in case \textnormal{(a)}, $\Pi$ is of type $(2,2,2)$
    if $\pi'$ is tetrahedral, and $(2, 4)$ otherwise.
\end{main}

\bigskip

\emph{Remark: } By \cite{KSh2000}, $\Pi = \pi \boxtimes \pi'$ is
automorphic on GL$(6)/F$, and its $L$--function agrees with the
Rankin--Selberg $L$--function $L (s, \pi \times \pi')$. Theorem A
implies in particular that $\Pi$ is cuspidal if (i) $\pi'$ is not
dihedral {\it and} (ii) $\pi$ is not a twist of $Ad(\pi')$.

\medskip

A partial cuspidality criterion was proved earlier by the second
author in \cite{Wa2003}; but he only treated the case when $\pi$
is twist equivalent to the Gelbart--Jacquet symmetric square
transfer of some cusp form on $GL (2)$.

\bigskip

\begin{main} \label{TM:B}
    Let $F$ be a number field and
    $\pi', \pi$ be cusp forms on $GL (2)/F$, $GL (3)/F$
    respectively. Put $\Pi = \pi \boxtimes \pi'$.
    Assume that $\pi'$ is not of solvable polyhedral type,
    and $\pi$ not essentially selfdual. Then we have the
    following:

    \textnormal{(a)} If $\pi$ does not
    admit any self twist, $\Pi$ is cuspidal
    without any self twist. Furthermore, if $\pi$
    is not monomial, then $\Pi$ is not induced from any
    non-normal cubic extension.

    \textnormal{(b)} If $\pi$ is not of solvable type,
    i.e., its base change to any solvable Galois extension
    is cuspidal, $\Pi$ is cuspidal and not of solvable
    type; in particular, there is no solvable extension
    $K / F$ of degree $d > 1$ dividing $6$, and a cuspidal automorphic
    representation $\eta$ of $GL_{6/d} (\A_{F})$,
    such that $L (s, \Pi) = L (s, \eta)$.
\end{main}

\medskip

\emph{Remark: } If $\pi$ is regular algebraic at infinity, and $F$
is not totally imaginary, then $\pi$ is not monomial (See Lemma
~\ref{T:903}).

\bigskip

We will prove Theorem ~\ref{TM:A} in sections 6 through 8, and
Theorem ~\ref{TM:B} in section 9.

Before proceeding with the proofs of these theorems, we will
digress and discuss the cohomological application.

\bigskip


\section{Preliminaries on cuspidal cohomology}

\bigskip

The experts can skip this section and go straight to the statement
(in section 4) and the proof (in section 5) of Theorems
~\ref{TM:C}, \ref{TM:D}. Let
$$
\Gamma \, \subset \, {\rm SL}(n, \Z),
$$
be a congruence subgroup of $G_n^0: = {\rm SL}(n, \R)$, which has
finite covolume. $\Gamma$ acts by left translation on the
symmetric space $X_n^0 : = \, {\rm SL}(n, \R)/{\rm SO}(n)$. The
cohomology of $\Gamma$ is the same as that of the locally
symmetric orbifold $\Gamma \backslash X_n^0$. If $H^\ast_{\rm
cont}$ denotes the {\it continuous group cohomology}, a version of
Shapiro's lemma gives an isomorphism
$$
H^\ast(\Gamma, \C) \, \simeq \, H^\ast_{\rm cont}(G_n^0, {\mathcal
C}^\infty(\Gamma\backslash G_n^0)).
$$
The constant functions are in this space, and the contribution of
$H^\ast_{\rm cont}(G_n^0, \C)$ to $H^\ast(\Gamma, \C)$ is well
understood and plays an important role in Borel's interpretation
of the values of the Riemann zeta function $\zeta(s)$ at negative
integers.

We will be interested here in another, more mysterious, piece of
$H^\ast(\Gamma, \C)$, namely its {\it cuspidal part}, denoted
$H^\ast_{\rm cusp}(\Gamma, \C)$, which injects into
$H^\ast(\Gamma, \C)$ by a theorem of Borel. Furthermore, one knows
by L.~Clozel (\cite{C}) that the cuspidal summand is defined over
$\Q$, preserved by the Hecke operators. The cuspidal cohomology is
represented by cocycles defined by smooth cusp forms in
$L^2(\Gamma\backslash G_n^0)$, i.e., one has
$$
H^\ast_{\rm cusp}(\Gamma, \C) \, = \, H^\ast_{\rm cont}(G_n^0,
L^2_{\rm cusp}(\Gamma\backslash G_n^0)^\infty),
$$
where $L^2_{\rm cusp}(\Gamma\backslash G_n^0)$ denotes the space
of cusp forms, and the superscript $\infty$ signifies taking the
subspace of smooth vectors. If $\mathfrak G_{n}$ denotes the {\it
complexified} Lie algebra of $G_n$, the passage from continuous
cohomology to the {\it relative Lie algebra cohomology}
(\cite{BoW}) furnishes an isomorphism
$$
H^\ast_{\rm cusp}(\Gamma, \C) \, \simeq \, H^\ast({\mathfrak
G}^0_{n},K:  L^2_{\rm cusp}(\Gamma\backslash G_n^0)^\infty).
$$
It is a standard fact (see \cite{BoJ}, for example) that the right
action of $G_n$ on $L^2_{\rm cusp}(\Gamma\backslash G_n^0)$ is
completely reducible, and so we may write
$$
L^2_{\rm cusp}(\Gamma\backslash G_n^0) \, \simeq \,
\hat{\oplus}_\pi \, m_\pi {\mathcal H}_\pi,
$$
where $\pi$ runs over the irreducible unitary representations of
$G_n^0$ (up to equivalence), ${\mathcal H}_\pi$ denotes the space
of $\pi$, $\hat{\oplus}$ signifies taking the Hilbert direct sum,
and $m_\pi$ is the multiplicity. Consequently,
$$
H^\ast_{\rm cusp}(\Gamma, \C) \, \simeq \, {\oplus}_\pi \,
H^\ast({\mathfrak G}_{n,\C}^0, K; {\mathcal
H}_\pi^\infty)^{m_\pi}.
$$
One knows completely which representations $\pi$ of $G_n^0$ have
non-zero $({\mathfrak G}_{n}^0, K)$-cohomology (\cite{VZ}; see
also \cite{Ku}). An immediate consequence (see \cite{C}, page 114)
is the following (with $[x]$ denoting, for any $x \in \R$, the
integral part of $x$):

\medskip

\begin{theorem} \label{T:301}
$$
H^i_{\rm cusp}(\Gamma, \C) \, = \, 0 \quad {\rm unless} \quad d(n)
\leq i \leq d(n)+[(n-1)/2],
$$
where
$$
d(n) = m^2 \quad {\rm if} \quad n=2m \quad {\rm and} \quad d(n) =
m(m+1) \quad {\rm if} \quad n=2m+1.
$$
\end{theorem}

\medskip

It will be necessary for us to work with the $\Q$-group $G_n: =
{\rm GL}(n)$ with center $Z_n$, and also work adelically. Let $\A
= \R \times \A_f$ be the adele ring of $\Q$, $K_\infty = O(n)$,
and $X_n = G_n(\R)/K_n$, whose connected component is $X_n^0$. For
any compact open subgroup $K$ of $G_n(\A_f)$, we have
$$
S_K: = \, G_n(\Q)Z_n(\R)^0\backslash G_n(\A)/K_\infty K \, \simeq
\, \cup_{j=1}^{r} \Gamma_j\backslash X_n^0, \leqno(3.6)
$$
where the $\Gamma_j$ are congruence subgroups of SL$(n, \Z)$ and
$Z_n(\R)^0$ is the Euclidean connected component of $Z_n(\R)$. We
need the following, which follows easily from the discussion in
section 3.5 of \cite{C}:

\begin{theorem} \label{T:302}
\begin{enumerate}
\item[(i)]
$$
H_{\rm cusp}^\ast(S_K, \C) \, \simeq \, \oplus_{\pi \in {\rm
Coh}_K} \, H^\ast(\tilde{\mathfrak G}_{n, \infty}, K_\infty;
\pi_\infty)\otimes \pi_f^K,
$$
where $\tilde{\mathfrak G}_{n, \infty}$ consists of matrices in
$M_n(\C)$ with purely imaginary trace, and {\rm Coh}$_K$ is the
set of (equivalence classes) of cuspidal automorphic
representations $\pi = \pi_\infty \otimes \pi_f$ of $G_n(\A)$ such
that $\pi_f^K \ne 0$, $\pi_\infty$ contributes to the relative Lie
algebra cohomology, and $(\omega_\pi)_\infty$ is trivial on
$Z(\R)^0$.

\item[(ii)] Suppose $\pi = \pi_\infty \otimes \pi_f$ is a cuspidal
automorphic representation of $G_n(\A)$ with $\pi_f^K \ne 0$ such
that the restriction $r_\infty$ of the Langlands parameter of
$\pi_\infty$ to $\C^\ast$ is given by the $n$-tuple
\begin{align}
\{ &(z/\vert z\vert)^{n-1}, ({\overline z}/\vert z\vert)^{n-1},
(z/\vert z\vert)^{n-3}, ({\overline z}/\vert z\vert)^{n-3},
 \notag \\
\dots, &(z/\vert z\vert), ({\overline z}/\vert z\vert)\} \otimes
(z\overline z)^{n-1}
\notag
\end{align}
if $n$ is even, and
\begin{align}
\{ &(z/\vert z\vert)^{n-1}, ({\overline z}/\vert z\vert)^{n-1},
(z/\vert z\vert)^{n-3}, ({\overline z}/\vert z\vert)^{n-3},
 \notag \\
\dots, &{(z/\vert z\vert)}^{2}, {({\overline z}/\vert z\vert)}^{2},
1\}
\otimes
(z\overline z)^{n-1}
\notag
\end{align}
if $n$ is odd. Then $\pi$ contributes to ${\rm Coh}_K$ in degree
$d(n)$.
\end{enumerate}
\end{theorem}

\qedsymbol

\bigskip

Given any cohomological $\pi$ as above, the fact that the cuspidal
cohomology at any level $K$ has a $\Q$-structure (\cite{C})
preserved by the action of the Hecke algebra ${\mathcal H}_\Q(G_f,
K)$ (consisting of $\Q$-linear combinations of $K$-double cosets),
implies that the $G_f$-module $\pi_f$ is {\it rational over a
number field} $\Q(\pi_f)$. When $n=2$, such a $\pi$ is defined by
a holomorphic newform $h$ of weight $2$, and then $\Q(\pi_f)$ is
none other than the field generated by the Fourier coefficients of
$h$.

\bigskip


\section{Non-selfdual, cuspidal classes for
$\Gamma \subset {\rm SL}(6, \Z)$}

\bigskip

The principle of functoriality predicts that given cuspidal
automorphic representations $\pi, \pi'$ of $G_n(\A), G_m(\A)$
respectively, there exists an isobaric automorphic representation
$\pi \boxtimes \pi'$ of $G_{nm}(\A)$ such that for every place $v$
of $\Q$, one has
$$
\sigma((\pi \boxtimes \pi')_v) \, \simeq \, \sigma(\pi_v) \otimes
\sigma_v(\pi'_v),
$$
where $\sigma$ is the map (up to isomorphism) given by the local
Langlands correspondence (\cite{HaT2000}, \cite{He2000})
from admissible irreducible
representations of $G_r(\Q_v)$ to $r$-dimensional representations
of $W'_v$, which is the real Weil group $W_\R$ if $v=\infty$ and
the extended Weil group $W_{\Q_p} \times {\rm SL}(2, \C)$ if $v$
is defined by a prime number $p$.

This prediction is known to be true for $n=m=2$ (\cite{Ra2000}),
More importantly for the matter at hand, it is also known for
$(n,m)=(3,2)$ by a difficult theorem of H.~Kim and F.~Shahidi
(\cite{KSh2000}).

\bigskip

Put
$$
T \, = \, T_1 \cup T_2,
$$
with
$$
T_1 = \{53, 61, 79, 89\} \quad {\rm and} \quad T_2 = \{128, 160,
205\}.
$$
By the article \cite{A-G-G} of Ash, Grayson and Green (for $p \in
T_1$), and the works \cite{vG-T1}, \cite{vG-T2}, \cite{vG-K-T-V}
of B.~van Geemen, J.~Top, et al (for $p \in T_2$), one knows that
for every $q \in T$, there is a non-selfdual cusp form $\pi(q)$ on
GL$(3)/\Q$ of level $q$, contributing to the (cuspidal) cohomology
(with constant coefficients) .

\bigskip

\begin{main} \label{TM:C}
    Let $\pi'$ be a cusp form on GL$(2)/\Q$ defined by a non-CM holomorphic newform
    $g$ of weight $4$, level $N$, trivial character, and
    field of coefficients $K$. Let $\pi$
    denote an arbitrary cusp form on GL$(3)/\Q$ contributing to
    the cuspidal cohomology in degree $2$, and let $\pi(q)$, $q \in T$,
    be one of the particular forms discussed above. Put $\Pi = \pi \boxtimes
    \pi'$ and $\Pi(q) = \pi(q) \boxtimes \pi'$. Then
    \begin{enumerate}
    \item[(a)] \, $\Pi$ contributes to the cuspidal cohomology of
    GL$(6)$.
    \item[(b)] \, $\Pi(q)$ is not essentially selfdual when $N \leq 23$ and $K = \Q$.
    \item[(c)] \, If $N$ is relatively prime to $q$, then the level of $\Pi(q)$ is
    $N^3q^2$. Now let $N \leq 23$ and $K=\Q$.
    Then $\Pi(q)$ does not admit any
    self-twist. Moreover, there is no cubic non-normal extension
    $K/\Q$ with a cusp form $\eta$ on GL$(2)/K$ such that $L(s, \Pi(q)) = L(s,
    \eta)$, nor is there a sextic extension (normal or not) $E/\Q$
    with a character $\lambda$ of $E$ such that $L(s, \Pi(q)) = L(s,
    \lambda)$.
    \end{enumerate}
\end{main}

\bigskip

We note from the {\it Modular forms database} of William Stein
(\cite{WSt}) that there exist newforms $g$ of weight $4$ with
$\Q$-coefficients, for instance for the levels $N = 5, 7, 13, 17,
19, 23$.

\medskip

\emph{Remark 4.1 } \, Part (b) should be true for any $\pi'$.
Suppose that for a given any cusp form $\pi$ on GL$(3)$,
cohomological or not, the functorial product $\Pi= \pi \boxtimes
\pi'$ satisfies $\Pi^\vee \simeq \Pi \otimes \nu$ for a character
$\nu$. Then at any prime $p$ where $a_p(\pi') \ne 0$, which
happens for a set of density $1$, we can of course conclude that
$a_p(\pi^\vee) = a_p(\pi)\nu(p)$. But this does not suffice, given
the state of knowledge right now concerning the refinement of the
strong multiplicity one theorem, to conclude that $\pi^\vee$ is
isomorphic to a twist of $\pi$. In the case of the $\pi(q)$, we
have information at a small set of primes and we have to make sure
that $a_p(\pi') \ne 0$ {\it and} $\overline{a_p(\pi)} \ne
a_p(\pi)$ for one of those $p$. The hypothesis that $K = \Q$ is
made for convenience, however, and the proof will extend to any
totally real field.

\bigskip

In \cite{vG-T1}, \cite{vG-T2}, \cite{vG-K-T-V} one finds in fact
an algebraic surface $S(q)$ over $\Q$ for each $q \in T_2$, and a
$3$-dimensional $\ell$-adic representation $\rho(q)$ (for any
prime $\ell$), occurring in $H^2_{\rm et}(S(q)_{\overline \Q},
\overline \Q_\ell)$, such that
$$
L_p(s, \rho(q)) \, = \, L_p(s, \pi(q)), \leqno(\ast)
$$
for all odd primes $p \leq 173$ not dividing $q$. Here is a
conditional result.

\bigskip

\begin{main} \label{TM:D}
    Let $\pi'$ be a cusp form on GL$(2)/\Q$ defined by a non-CM holomorphic newform
    $g$ of weight $4$, level $N$, trivial character, and field of coefficients $\Q$,
    with corresponding $\overline \Q_\ell$-representation $\rho'$.
    Let $\pi(q), \rho(q), \Pi(q)$ be as above for $q \in
    T_2$. Put $R(q) = \rho(q) \otimes \rho'$. Suppose $(\ast)$ holds at all the odd primes
    $p$ not dividing $q\ell$. Then $R(q)$ remains irreducible when
    restricted to any open subgroup of Gal$(\overline \Q/\Q)$.
\end{main}

\bigskip


\section{Proof of Theorems C, D modulo Theorems A, B}

\bigskip

In this section we will assume the truth of Theorems ~\ref{TM:A} and ~\ref{TM:B}.

{\it Proof of Theorem ~\ref{TM:C}}: \, As $\pi'$ is associated to
a holomorphic newform of weight $4$, we have
$$
\sigma(\pi'_\infty)\vert_{\C^\ast} \, \simeq \, \left((z/|z|)^3
\oplus (\overline z/|z|)^3\right) \otimes (z\overline z)^3.
$$
And since $\pi$ contributes to cohomology, we have (cf. part (ii)
of Theorem ~\ref{T:301})
$$
\sigma(\pi_\infty)\vert_{\C^\ast} \, \simeq \, \left((z/|z|)^2
\oplus 1 \oplus (\overline z/|z|)^2\right) \otimes (z\overline
z)^2.
$$
Since $\Pi_\infty$ corresponds to the tensor product
$\sigma(\pi_\infty) \otimes \sigma(\pi'_\infty)$, we get part (a)
of Theorem ~\ref{TM:C} in view of Theorem ~\ref{TM:A} and part
(ii) of Theorem ~\ref{T:302}.

\bigskip

Pick any $q$ in $T$ and denote by $\Q(\pi(q))$ the field of
rationality of the finite part $\pi(q)_f$ of $\pi(q)$. Then it is
known by \cite{A-G-G} that for $q \in T_1$,
$$
\Q(\pi(53))=\Q(\sqrt{-11}), \, \Q(\pi(61)) = \Q(\sqrt{-3}), \,
\Q(\pi(79)) = \Q(\sqrt{-15}), \, \Q(\pi(89)) = \Q(i),
$$
while by \cite{vG-T1}, \cite{vG-K-T-V} and \cite{vG-T2},
$$
\Q(\pi(q)) \, = \, \Q(i), \, \, \, \forall \, \, q \in T_2.
$$

By hypothesis, $\pi'$ is non-CM, and by part (a), $\Pi(q)$ is
cuspidal. Suppose there exists a character $\nu$ such that for
some $q \in T$,
$$
\Pi(q)^\vee \, \simeq \, \Pi(q) \otimes \nu.
$$
Comparing central characters, we get $\nu^6 = 1$. We claim that
$\nu^2 = 1$. Suppose not. Then there exists an element $\sigma$ of
Gal$(\overline \Q/\Q)$ fixing $\Q(\pi(q))$ such that $\nu \ne
\nu^\sigma$. Since $\pi'$ has $\Q$-coefficients and $\pi(q)$ has
coefficients in $\Q(\pi(q))$, we see that $\Pi(q)_{f}$ must be
isomorphic to the Galois conjugate $\Pi(q)_{f}^\sigma$, which
exists because the cuspidal cohomology group has, by Clozel (see
section 3), a $\Q$-structure preserved by the Hecke operators. If
we put $\mu = \nu/\nu^\sigma \ne 1$, we then see that $\Pi(q)
\simeq \Pi(q) \otimes \mu$. But we will see below that $\Pi(q)$
admits no non-trivial self-twist. This gives the desired
contradiction, proving the claim. If $\nu$ is non-trivial, the
quadratic extension $F/\Q$ it cuts out will need to have
discriminant dividing $q^aN^b$ for suitable integers $a, b$. For
any prime $p$ which is unramified in $F$, we will have
$$
\overline a_p(\pi)a_p(\pi') \, = \, \pm a_p(\pi)a_p(\pi').
$$
For each $j \leq 3$ and for each $\pi'$ with $N \leq 23$ and $K =
\Q$, we can find, using the tables in \cite{A-G-G},
\cite{vG-T1},\cite{vG-T2}, \cite{vG-K-T-V} and \cite{WSt}, a prime
$p$ such that $a_p(\pi') \ne 0$, $\nu(p) \ne 0$ and $\overline
a_p(\pi) \ne a_p(\pi)$. This proves part (b) of Theorem
~\ref{TM:C}.

\bigskip

When $N$ is relatively prime to $q$, the conductor of $\Pi(q)$
must be $N^3q^2$ as can be seen from the way epsilon factors
change under twisting (see section 4 of \cite{Ba-R} for example).

{F}rom now on, let $N \leq 23$ and $K = \Q$. One knows that as
$\pi'$ is holomorphic and not dihedral, the associated Galois
representation $\rho'$ remains irreducible when restricted to any
open subgroup of Gal$(\overline \Q/\Q)$. It follows that the base
change of $\pi'$ to any solvable Galois extension remains
cuspidal. In particular, it is not of solvable polyhedral type. We
claim that $\pi(q)$ is not monomial. Indeed, the infinite type of
$\pi(q)$ is regular algebraic \cite{C}, and to be monomial there
needs to be a cubic, possibly non-normal, extension $K/\Q$ which
can support an algebraic Hecke character which is {\it not} a
finite order character times a power of the norm. By \cite{Weil},
for such a character to exist, $K$ must contain a CM field, i.e.,
a totally imaginary quadratic extension of a totally real field,
which forces $K$ to be imaginary. But any cubic extension of $\Q$
has a real embedding, and this proves our claim. Note also that as
$\pi(q)$ is not essentially self-dual, it is not a twist of the
symmetric square of any cusp form, in particular $\pi'$, on
GL$(2)/\Q$. Now it follows from Theorem ~\ref{TM:B} that $\Pi(q)$
does not admit any self-twist.

Suppose $K$ is a non-normal cubic field together with a cusp form
$\eta$ on GL$(2)/K$ such that $L(s, \Pi(q)) = L(s, \eta)$. Let $L$
be the Galois closure of $K$ (with Galois group $S_3$), and let
$E$ be the quadratic extension of $\Q$ contained in $L$. Then
$\Pi(q)_{E}$ will be cuspidal and automorphically induced by the
cusp form $\eta_L$ of GL$(3, \A_L)$. In other words, $\Pi(q)_{E}$
admits a non-trivial self-twist. To contradict this, it suffices,
in view of Theorem ~\ref{TM:B}, to show that $\pi(q)_{E}$ admits
no self-twist relative to $L/E$, i.e., that $\pi(q)_{E}$ is not
automorphically induced by a character $\mu$ of $L$. But as noted
above, this forces $L$ to be a totally imaginary number field
containing a CM field $L_0$. Then either $L = L_0$ or $L_0 = E$.
In the latter case, by \cite{Weil}, $\mu$ will be a finite order
character times the pullback by norm of a character $\mu_0$ of
$E$, forcing $I_L^E(\mu)$ to be {\it not} regular at infinity, and
so this case cannot happen. So $L$ itself must be a CM field, with
its totally real subfield $F$. Then Gal$(F/\Q)$ would be cyclic of
order $3$ and a quotient of $S_3$, which is impossible. So this
case does not arise either. So $\pi(q)_{E}$ doe snot admit any
self-twist, and $\Pi(q)$ is not associated to any $\eta$ as above.

Now suppose $L(s, \Pi(q)) = L(s, \lambda)$ for a character
$\lambda$ of a sextic field $L$. If $L$ contains a proper subfield
$M \ne \Q$, then since $m:= [L:M] \leq 3$, one can induce
$\lambda$ to $M$ and get an automorphic representation $\beta$ of
GL$_m(\A_{M})$ such that $L(s, \lambda) = L(s, \beta) = L(s,
\Pi(q))$, which is impossible by what we have seen above. So $L$
must not contain any such $M$. But on the other hand, since
$\Pi(q)_{\infty}$ is algebraic and regular, we need $L$ to
contain, by \cite{Weil}, a CM subfield $L_0$, and hence also its
totally real subfield $F$. Either $F = \Q$, in which case $L_0$ is
imaginary quadratic, or $F \ne \Q$. Either way there will be a
proper subfield $M$ of degree $\leq 3$, and so the purported
equality $L(s, \Pi(q)) = L(s, \lambda)$ cannot happen. We are now
done with the proof of Theorem ~\ref{TM:C}.

\bigskip

{\it Proof of Theorem ~\ref{TM:D}}: \, By assumption, the
$\ell$-adic representation $\rho$ is functorially associated to
the cuspidal cohomological form $\pi(q)$ on GL$(3)/\Q$ with $q \in
T_2$.

\medskip

\begin{lemma} \label{T:501}
$\rho$ is irreducible under restriction to any open subgroup.
\end{lemma}

\medskip

{\it Proof}. \, Suffices to show that the restriction $\rho_E$ to
Gal$(\overline \Q/E)$ is irreducible for any finite {\it Galois}
extension $E/\Q$. Pick any such extension and write $G =
$Gal$(E/\Q)$. Suppose $\rho_E$ is reducible. Then we have {\it
either}
\begin{enumerate}
\item[(i)] \, $\rho_E \simeq \tau \oplus \chi$ with $\tau$ irreducible of
dimension $2$ and $\chi$ of dimension $1$; {\it or}
\item[(ii)] \, $\rho_E \simeq \chi_1 \oplus \chi_2 \oplus \chi_3$,
with each $\chi_j$ one-dimensional.
\end{enumerate}
Let $V$ be the $3$-dimensional $\overline \Q_\ell$-vector space on
which Gal$(\overline \Q/\Q)$ acts via $\rho$. Suppose we are in
case(i), so that there is a line $L$ in $V$ preserved by
Gal$(\overline \Q/E)$ and acted upon by $\chi$. Note that $G$ acts
on $\{\tau, \chi\}$ and, by the dimension consideration, it must
preserve $\{\chi\}$. Hence the line $L$ is preserved by
Gal$(\overline \Q/\Q)$, which contradicts the fact that $\rho$ is
irreducible.

So we may assume that we are in case (ii). We claim that $\chi_i
\ne \chi_j$ if $i \ne j$. Indeed, since $\rho$ arises as (the base
change to $\overline \Q_\ell$ of) a summand of the $\ell$-adic
cohomology of a smooth projective variety, it is Hodge-Tate, and
so is each $\chi_j$. So each $\chi_j$ is locally algebraic and
corresponds to an algebraic Hecke character $\chi_j'$ of $E$. By
the identity of the $L$-functions, we will have $L^S(s, \pi) =
\prod_j L^S(s, \chi'_j)$ for a suitable finite set $S$ of places
$S$. By the regularity of $\pi$, each $\chi_j'$ must appear with
multiplicity one, which proves the claim. Now let $L_j$ denote,
for each $j \leq 3$, the (unique) line in $V$ stable under
Gal$(\overline \Q/E)$ and acted upon by $\chi_i$. And $G$ acts by
permutations on the set $\{\chi_1, \chi_2, \chi_3\}$. In other
words, there is a representation $r: G \rightarrow S_3$ such that
the $G$-action is via $r$. Put $H = {\rm Ker}(r)$, with
corresponding intermediate field $M$. Then each $L_j$ is stable
under Gal$(\overline \Q/M)$, so that $\rho_M \simeq \nu_1 \oplus
\nu_2 \oplus \nu_3$, where each $\nu_j$ is a character of
Gal$(\overline \Q/M)$. Also, $M/\Q$ is Galois with Gal$(M/\Q)
\subset S_3$. But from the proof of Theorem ~\ref{TM:C} that the
base change $\pi_M$  of $\pi$ to any such $M$ is cuspidal.
However, if $\nu_j^1$ denotes the algebraic Hecke character of $M$
defined by $\nu_j$, the twisted $L$-function $L^S(s, \pi_M \otimes
{\nu'_j}^{-1})$ will have a pole at $s=1$, leading to a
contradiction. We have now proved Lemma 5.1.

\qed

\medskip

Note that Lemma 1 implies in particular that for any finite
extension $F/\Q$, $\rho_F$ does not admit any self-twist.

\bigskip

\begin{lemma} \label{T:502}
For any finite extension $E/\Q$, the restriction $\rho_E$ is not
essentially self-dual.
\end{lemma}

\medskip

{\it Proof}. \, Again we may assume that $E/\Q$ is Galois with
group $G$. As before let $V$ denote the space of $\rho$, and
suppose that we have an isomorphism $\rho \simeq \rho^\vee \otimes
\nu$, for a character $\nu$. Then there is a line $L$ in $V
\otimes V$ on which Gal$(\overline \Q/E)$ acts via $\nu$. By
Schur's lemma (and this is why we have to work over $\overline
\Q_\ell$), the trivial representation appears with multiplicity
one in $V \otimes V^\vee$. It implies that $\nu$ must appear with
multiplicity one in $V \otimes V$. We claim that $V \otimes V$
contains no other character. Indeed, if we have another character
$\nu'$, we would have $\rho \simeq \rho \otimes \mu$, where $\mu =
\nu/\nu'$. But as noted above, $\rho_E$ admits no self-twist, and
so $\mu = 1$, and the claim is proved. Consequently, the action of
$G$ on $V \otimes V$ must preserve $\nu$. In other words, the line
$L$ is stable under all of Gal$(\overline \Q/\Q)$, contradicting
the fact that $\rho$ is not essentially self-dual. Done

\qed

Now consider $R = \rho \otimes \rho'$. We know that both $\rho$
and $\rho'$ remain irreducible upon restriction to any open
subgroup and moreover, such a restriction of $\rho$ is {\it not}
essentially self-dual. It then follows easily that the restriction
of $R$ is irreducible.

This finishes the proof of Theorem ~\ref{TM:D}.

\qed

\bigskip


\section{Proof of Theorem A, Part \#1}

By twisting we may assume that $\pi, \pi'$ are unitary, so that
$\pi^\vee \simeq \overline \pi$ and ${\pi'}^\vee \simeq
\overline{\pi'}$, with respective central characters $\omega,
\omega'$.

Now we proceed in several steps. Applying Langlands's
classification, (\cite{La79-1}, \cite{La79-2}, \cite{JS81}), we
see that the Kim-Shahidi product $\Pi= \pi \boxtimes \pi'$ must be
an isobaric sum of cusp forms whose degrees add up to $6$. Thanks
to the Clebsch-Gordon decomposition
$$
{\rm sym}^2(\pi') \boxtimes \pi' \, \simeq \, {\rm sym}^3(\pi')
\boxplus (\pi' \otimes \omega'),
$$
$\Pi$ is not cuspidal if $\pi$ is a twist of ${\rm sym}^{2}
(\pi')$.

The list of all the cases when $\Pi$ is not cuspidal is the
following:

\medskip

{\bf Case I}: \, {\it $\Pi$ has a constituent of degree $1$, i.e.,
$\Pi = \lambda \boxplus \Pi'$ for some idele class character
$\lambda$ and some automorphic representation $\Pi'$ of $GL (5)$.}

\medskip

{\bf Case II}: \, {\it $\Pi$ has a constituent of degree $2$,
i.e., $\Pi = \tau \boxplus \Pi'$ for some cusp form $\tau$ on $GL
(2)$ and some automorphic representation $\Pi'$ of $GL (4)$.}

\medskip

{\bf Case III}: \, {\it $\Pi$ is an isobaric sum of two cusp forms
$\sigma_{1}$ and $\sigma_{2}$ on $GL (3)$. }

\bigskip

We first deal with Cases I and II. We need some preliminaries.
First comes the following basic result due to H. ~Jacquet and J.A.
~Shalika (\cite{JS81}, \cite{JS90}, and R. ~Langlands
(\cite{La79-1}, \cite{La79-2}).

\bigskip

\begin{lemma} \label{T:601}

\textnormal{(i)} Let $\Pi$, $\tau$ be isobaric automorphic
representations of $GL_{n} (\mathbb{A}_{F})$, $GL_{m}
(\mathbb{A}_{F})$ respectively. Assume that $\tau$ is cuspidal.
Then the order of the pole of $L (s, \Pi \otimes \overline{\tau})$
at $s=1$ is the same as the multiplicity of $\tau$ occurring in
the isobaric sum decomposition of $\Pi$.

\textnormal{(ii)} $L (s, \Pi \times \overline{\Pi})$ has a pole at
$s = 1$ of order $m = \sum_{i} m_{i}^{2}$ if $\Pi = \boxplus_{i}
m_{i} \pi_{i}$ is the isobaric decomposition of $\Pi$, with the
$\pi_{i}$ being inequivalent cuspidal representations of smaller
degree. In particular, $m = 1$ if and only if $\Pi$ is cuspidal.
\end{lemma}

\bigskip

An $L$--function $L(s)$ is said to be nice if it converges on some
right half plane, admits an Euler product of some degree $m$, say,
and extends to a meromorphic function of finite order with no pole
outside $s = 1$, together with a functional equation related to
another $L$--function  $L^{\vee} (s)$ given by
\[
    L (s) = W {(d_{F}^{m} N)}^{1/2 - s} L^{\vee} (1 - s),
\]
where $W$ is a non-zero scalar.

\medskip

If $\pi_{1}, \pi_{2}$ are automorphic forms on $GL (m), GL (n)$
respectively, then the Rankin-Selberg $L$-function $L (s, \pi_{1}
\times \pi_{2})$ is known to be nice (\cite{JPSS}, \cite{MW},
\cite{Sh}). Of course, the product of two nice $L$--functions is
nice. Furthermore, we recall the following Tchebotarev-like result
for nice $L$--functions (\cite{JS90}):

\begin{lemma} \label{T:602}
    Let $L_{1} (s) = \prod_{v} L_{1, v} (s)$ and $L_{2} (s) = \prod_{v} L_{2, v} (s)$
    be two $L$--functions with Euler products, and suppose that they are both of exactly one of the following types:

    \textnormal{(a)} $L_{i} (s)$ is an Artin $L$--function of some Galois extension;

    \textnormal{(b)} $L_{i} (s)$ is attached to an isobaric automorphic representation;

    \textnormal{(c)} $L_{i} (s)$ is a Rankin--Selberg $L$--function of two isobaric
    automorphic representations.

    If $L_{1, v} (s) = L_{2, v} (s)$ for all but finite places $v$ of
    $F$, then $L_{1} (s) = L_{2} (s)$.
\end{lemma}

\bigskip

\emph{Proof of Theorem A for Cases I and II.}

\medskip

Firstly, Case I can never happen. The reason is the following: If
$\lambda$ is a constituent of $\Pi = \pi \boxtimes \pi'$, then $L
(s, \Pi \times \bar{\lambda})$ has a pole at $s = 1$ (Lemma
~\ref{T:601}), hence so does
$$
L (s, \pi' \times \pi \otimes \bar{\lambda}) = L (s, \pi \boxtimes
\pi' \otimes \bar{\lambda}).
$$
However, $\pi'$ and $\pi \otimes \bar{\lambda}$ are cuspidal of
different degrees, hence $L (s, \pi' \times \pi \otimes
\bar{\lambda})$ is entire, and we get the desired contradiction.

\medskip

Now we treat Case II, where $\Pi$ has a constituent $\tau$ of
degree $2$. We will show that this can happen IF AND ONLY IF $\pi$
is twist equivalent to ${\rm sym}^{2} (\pi')$ in which case $\tau$
is twist equivalent to $\pi'$.

\medskip

In fact, for each finite $v$ where $\pi$ and $\pi'$ are
unramified,
$$
L_{v} (s, \Pi \otimes \bar{\tau}) = L_{v} (s, \pi
\times (\pi' \boxtimes \bar{\tau})),
$$
where $\pi' \boxtimes \bar{\tau}$ is the functorial product of
$\pi'$ and $\bar{\tau}$ whose modularity (in GL$(4)$) was
established in \cite{Ra2000}. One may check the following: If
$\pi'_{v} = \alpha_{v, 1} \boxplus \alpha_{v, 2}$, $\pi_{v} =
\beta_{v, 1} \boxplus \beta_{v, 2} \boxplus \beta_{v, 3}$, and
$\tau'_{v} = \gamma_{v, 1} \boxplus \gamma_{v, 2}$, then both
sides of the equality is the same as $\prod_{i, j, k} L (s,
\alpha_{v, i} \beta_{v, j} \bar{\gamma}_{v, k})$ where the product
is over all $i, j, k$ such that $1 \leq i, k \leq 2$ and $1 \leq j
\leq 3$.

Hence by Lemma ~\ref{T:602},
\[
    L (s, \Pi \otimes \bar{\tau}) = L (s, \pi \times (\pi' \boxtimes \bar{\tau}))
\]

\medskip

As $\tau$ is a constituent of $\Pi$, the $L$--functions on both
sides above have a pole at $s = 1$. As $\pi$ is cuspidal, this
means by Lemma ~\ref{T:601}, $\bar{\pi}$ is a constituent of $\pi'
\boxtimes \bar{\tau}$. Hence $\pi' \boxtimes \bar{\tau}$ should
possess a constituent of degree $1$, namely a character $\mu$.

\medskip

Thus $L (s, \pi' \times \bar{\tau} \otimes \bar{\mu}) = L (s, \pi'
\boxtimes \bar{\tau} \otimes \bar{\mu})$ has a pole at $s = 1$,
implying that $\pi'$ is equivalent to $\tau \otimes \mu$. Hence
$$
\pi' \boxtimes \bar{\tau} \cong \mu \boxplus Ad (\tau) \otimes
\mu,
$$
which means that $\pi \cong Ad (\tau) \otimes \mu \cong Ad
(\pi') \otimes \mu$.

\medskip

Finally, it is clear that if Case II happens, then $\pi'$ cannot
be dihedral. Furthermore, $\Pi$ is Eisensteinian of type $(2, 2,
2)$ if $\pi'$ is tetrahedral, and $(2, 4)$ otherwise. we can see
this by observing that
\[
    \pi' \boxtimes Ad (\pi') \cong {\rm sym}^{3} (\pi') \otimes
\omega_{\pi'} \boxplus \pi' \otimes \omega_{\pi'}^{2}
\]

\qedsymbol

\bigskip


\section{Proof of Theorem A, Part \#2}

It remains to treat Case III. Here again, $\pi'$ denotes a cusp
form on $GL (2)$ and $\pi$ a cusp form on $GL (3)$. Assume that
$\Pi = \sigma_{1} \boxplus \sigma_{2}$ where $\sigma_{1}$ and
$\sigma_{2}$ are cusp forms on $GL (3)$.

\medskip

We will divide Case III into two subcases: In this section, we
will assume that $\pi'$ is not dihedral. The (sub)case when $\pi'$
is dihedral will be treated in the next section.

\bigskip

The following equality is crucial, and it holds for all cusp forms
$\pi'$ on $GL (2)$ and $\pi$ on $GL (3)$:

\medskip

\begin{proposition}  \label{T:701}
\begin{align}
    L (s, &\pi \times \pi'; \Lambda^{3} \otimes \omega_{\pi}^{-1} \chi)
      L (s, \pi' \otimes \omega_{\pi'} \chi)
    \notag \\ \label{EQ:701}
    &= L (s, {\rm sym}^{3} (\pi') \otimes \chi)
     L (s, (\pi \boxtimes \pi') \times \tilde{\pi} \otimes \omega_{\pi'} \chi)
\end{align}
where $\omega_{\pi'}$ and $\omega_{\pi}$ are the respective
central characters of $\pi'$ and $\pi$.
\end{proposition}

\medskip

\emph{Proof of Proposition ~\ref{T:701}.}

We claim that both sides of \eqref{EQ:701} are nice. Indeed, we
see that formally, the admissible representation $\Lambda^{3} (\pi
\boxtimes \pi') \otimes \omega_{\pi}^{-1}$ is equivalent to ${\rm
sym}^{3} (\pi') \boxplus (Ad (\pi) \boxtimes \pi' \otimes
\omega_{\pi'})$. So the left hand side is nice. And the right hand
side is nice by \cite{KSh2000}, whence the claim. So by Lemma
~\ref{T:602} it suffices to prove this equality given by the
Proposition locally at $v$ for almost all $v$. It then suffices to
prove the following identity (as admissible representations) for
almost all $v$:
\begin{align}
    \Lambda^{3} &(\pi'_{v} \boxtimes \pi_{v}) \otimes \omega_{\pi_{v}}^{-1}
    \boxplus \pi'_{v} \otimes \omega_{\pi'_{v}}
    \notag \\ \label{EQ:702}
    &= {\rm sym}^{3} (\pi'_{v}) \boxplus \pi'_{v} \boxtimes
    \pi_{v} \boxtimes \tilde{\pi}_{v} \otimes \omega_{\pi'_{v}}
\end{align}

\medskip

Let $v$ be any place where $\pi'$ and $\pi$ are unramified. Say
$\pi'_{v} = \alpha_{v, 1} \boxplus \alpha_{v, 2}$, $\pi_{v} =
\beta_{v, 1} \boxplus \beta_{v, 2} \boxplus \beta_{v, 3}$. Note
that $\omega_{\pi'_{v}} = \alpha_{v, 1} \alpha_{v, 2}$ and
$\omega_{\pi_{v}} = \beta_{v, 1} \beta_{v, 2} \beta_{v, 3}$. Then
it is routine to check that the left and the right hand sides of
\eqref{EQ:702} are equal to the sum of the following terms:

\medskip

Terms A ($2$ terms): \, $\alpha_{v, 1}^{3} \boxplus \alpha_{v,
2}^{3}$;

Terms B ($2 \times 4 = 8$ terms): \, $4$ copies of $(\alpha_{v, 1}
\boxplus \alpha_{v, 2}) \otimes \omega_{\pi'_{v}}$;

Terms C ($12$ terms): \, $\boxplus_{1 \leq i \leq 2, 1 \leq j \ne
k \leq 3}\, \alpha_{v, i} \omega_{\pi'_{v}} \beta_{v, j}
\beta^{-1}_{v, k}$.

\medskip

In fact, Terms A, B and C are obtained by expanding the right hand
side of \eqref{EQ:702}. Since
\[
    Ad (\pi_{v}) = 3 \cdot \underline{1} \boxplus (\boxplus_{1 \leq j \ne k \leq 3},
    \beta_{v, j} \beta^{-1}_{v, k}),
\]
the Terms C and (three of) the Terms B are obtained from $\pi'_{v}
\boxtimes \pi_{v} \boxtimes \tilde{\pi}_{v} \otimes
\omega_{\pi'}$, and the Terms A and (one of) the Terms B arise
from ${\rm sym}^{3} (\pi')$.

\medskip

The left hand side is easy to handle since we have the following:

\begin{align}
    \Lambda^{3} (\pi'_{v} \otimes \pi_{v}) &=
    \boxplus_{1 \leq i \leq 2, 1 \leq j, k \leq 3, j \ne k}\,
    \alpha_{v, i} \omega_{\pi'_{v}} \beta_{v, j}^{2} \beta_{v, k}
    \notag \\
    &\boxplus \omega_{\pi_{v}} \alpha_{v, 1}^{3}
    \boxplus \omega_{\pi_{v}} \alpha_{v, 2}^{3}
    \notag \\
    &\boxplus 3 \omega_{\pi_{v}} \alpha_{v, 1} \omega_{\pi'_{v}}
    \boxplus 3 \omega_{\pi_{v}} \alpha_{v, 2} \omega_{\pi'_{v}}
    \notag
\end{align}

In fact, the $\omega_{\pi_{v}}^{-1}$ twist of the thing above
contributes the Terms C, Terms A and (three of) Terms B.

\medskip

So we have proved \eqref{EQ:702}, and hence \eqref{EQ:701}.

\qedsymbol

\bigskip

Let $\sigma_{1}$ and $\sigma_{2}$ be cusp forms on $GL (3)$.

\begin{lemma} \label{T:702}
    Let $\eta_{1}$ and $\eta_{2}$ be the central characters of $\sigma_{1}$
    and $\sigma_{2}$ respectively.
    Then
\begin{align}
    L (s, \sigma_{1} \boxplus \sigma_{2}; \Lambda^{3} \otimes \chi')
    &= L (s, \eta_{1} \chi') L (s, \eta_{2} \chi')
    \notag \\ \label{EQ:703}
    &L (s, \sigma_{1} \times \tilde{\sigma}_{2} \otimes \eta_{2} \chi')
    L (s, \sigma_{2} \times \tilde{\sigma}_{1} \otimes \eta_{1} \chi')
\end{align}
\end{lemma}

\emph{Proof. of Lemma ~\ref{T:702}.}

This is easy since at each place $v$ where the $\sigma_{i}$ are
unramified,
\[
    \Lambda^{3} (\sigma_{1, v} \boxplus \sigma_{2, v})
    = \boxplus_{0 \leq i \leq 3}\, \left( \Lambda^{i} (\sigma_{1, v}) \boxtimes
    \Lambda^{3 - i} (\sigma_{2, v}) \right),
\]

\[
    \Lambda^{2} (\sigma_{i, v}) \cong \tilde{\sigma}_{i, v} \otimes \eta_{i}
\]
and
\[
    \Lambda^{3} (\sigma_{i}) \cong \eta_{i}.
\]
Done by applying Lemma ~\ref{T:602}.

\qedsymbol

\bigskip

Before we apply Proposition ~\ref{T:701} and Lemma ~\ref{T:702},
let us first investigate a special instance of Case III when
$\sigma_{1}$ and $\sigma_{2}$ are both twists of $\pi$:

\begin{lemma} \label{T:703}
    If $\pi \boxtimes \pi' = (\pi \otimes \chi_{1}) \boxplus (\pi \boxplus \chi_{2})$
    then
\[
    {\rm sym}^{3} (\pi') \cong (\pi' \otimes \omega_{\pi'}) \boxplus \chi_{1}^{3}
    \boxplus \chi_{2}^{3}
\]
    Hence if $\pi'$ is not dihedral or tetrahedral, this cannot happen.
\end{lemma}

\emph{Proof of Lemma ~\ref{T:703}.}

Let $v$ be any place where $\pi'_{v}$ and $\pi_{v}$ are
unramified. Write $$L_{v} (s, \pi') = {(1 - U_{v}
{(Nv)}^{-s})}^{-1} {(1 - V_{v} {(Nv)}^{-s})}^{-1}$$ and $$L_{v}
(s, \pi) = {(1 - A_{v} {(Nv)}^{-s})}^{-1} {(1 - B_{v}
{(Nv)}^{-s})}^{-1} {(1 - C_{v} {(Nv)}^{-s})}^{-1}.$$

\medskip

Then
\begin{align}
    L_{v} &(s, \pi' \times \pi) =
    \notag \\
    &     {(1 - U_{v} A_{v} {(Nv)}^{-s})}^{-1} {(1 - U_{v}  B_{v} {(Nv)}^{-s})}^{-1}
    {(1 - U_{v} C_{v} {(Nv)}^{-s})}^{-1}
    \notag \\
    &     {(1 - V_{v} A_{v} {(Nv)}^{-s})}^{-1} {(1 - V_{v}  B_{v} {(Nv)}^{-s})}^{-1}
    {(1 - V_{v} C_{v} {(Nv)}^{-s})}^{-1}
    \notag
\end{align}

\medskip

Let $X_{v} = \chi_{1} ({\rm Frob}_{v})$ and $Y_{v} = \chi_{2}
({\rm Frob}_{v})$. Then
\begin{align}
    L_{v} &(s, \pi \otimes \chi_{1}) L_{v} (s, \pi \otimes \chi_{2}) =
    \notag \\
    &     {(1 - X_{v} A_{v} {(Nv)}^{-s})}^{-1} {(1 - X_{v}  B_{v} {(Nv)}^{-s})}^{-1}
    {(1 - X_{v} C_{v} {(Nv)}^{-s})}^{-1}
    \notag \\
    &     {(1 - Y_{v} A_{v} {(Nv)}^{-s})}^{-1} {(1 - Y_{v}  B_{v} {(Nv)}^{-s})}^{-1}
    {(1 - Y_{v} C_{v} {(Nv)}^{-s})}^{-1}
    \notag
\end{align}

\medskip

Consequently we have
\begin{align}
    &\Set{U_{v} A_{v}, U_{v} B_{v}, U_{v} C_{v}, V_{v} A_{v}, V_{v} B_{v}, V_{v} C_{v}} =
   \notag \\ \label{EQ:704}
    &\Set{X_{v} A_{v}, X_{v} B_{v}, X_{v} C_{v}, Y_{v} A_{v}, Y_{v} B_{v}, Y_{v}
    C_{v}}.
\end{align}

\medskip

In particular,
\[
(U_{v}^{n} + V_{v}^{n})(A_{v}^{n} + B_{v}^{n} + C_{v}^{n})
=(X_{v}^{n} + Y_{v}^{n})(A_{v}^{n} + B_{v}^{n} + C_{v}^{n})
\]
for each positive integer $n$. And besides, taking the products of
the elements on each side of \eqref{EQ:704} and equating, we get
\[
    U_{v}^{3} V_{v}^{3} = X_{v}^{3} Y_{v}^{3}
\]

\medskip

Now we apply the following lemma.

\bigskip

\begin{lemma} \label{T:707}
    If $X, Y, U, V, A, B, C$ are nonzero complex numbers
    such that (for all $n > 0$)
    \[
        (U^{n} + V^{n})(A^{n} + B^{n} + C^{n})
        = (X^{n} + Y^{n})(A^{n} + B^{n} + C^{n}),
    \]
    and
    $U^{3} V^{3} = X^{3} Y^{3}$,
    then $\Set{U^{3}, V^{3}} = \Set{X^{3}, Y^{3}}$.
\end{lemma}

\emph{Proof of Lemma ~\ref{T:707}.}

If $A^{3} + B^{3} + C^{3} \ne 0$, then $U^{3} + V^{3} = X^{3} +
Y^{3}$. Hence $\Set{U^{3}, V^{3}} = \Set{X^{3}, Y^{3}}$ as $U^{3}
V^{3} = X^{3} Y^{3}$.

\medskip

If $A^{3} + B^{3} + C^{3} = 0$, we claim that $A + B + C \ne 0$.
Otherwise,
$$- 3 A B C = 3 A B (A + B) = {(A + B)}^{3} - (A^{3} + B^{3}) = -C^{3} + C^{3} = 0$$
Thus $a$, $b$ or $c$ is zero. This leads to a contradiction.

\medskip

In fact we will prove the following statement:

{\bf Claim:} \emph{If $a, b, c \ne 0$, then $a + b + c$ or $a^{3}
+ b^{3} + c^{3}$ is not zero.}

\medskip

So we claim also that $A^{9} + B^{9} + C^{9} \ne 0$, and, $A^{2} +
B^{2} + C^{2}$ or $A^{6} + B^{6} + C^{6}$ is not zero.

\medskip

Hence
\[
    U^{n} + V^{n} = X^{n} + Y^{n}
\]
holds for $n = 1$ and $9$, and for one of $2$ or $6$.

\medskip

If this equality holds for $n = 1$, and $2$, then
\[
    U V = \frac{{(U + V)}^{2} - (U^{} + V^{2})}{2}
    = \frac{{(X + Y)}^{2} - (X^{} + Y^{2})}{2} = X Y,
\]
implying that $\Set{U, V} = \Set{X, Y}$, and the lemma will
follow.

\medskip

Now assume that $U^{n} + V^{n} = X^{n} + Y^{n}$ holds for $n = 1$,
$6$ or $9$. As we have already assumed that $U^{3} V^{3} = X^{3}
Y^{3}$, $U^{3 n} V^{3 n} = X^{3 n} Y^{3 n}$. So we have
\[
    \Set{U^{n}, V^{n}} = \Set{X^{n}, Y^{n}}
\]
for $n = 6$ and $9$.

\medskip

Without loss of generality, assume that $U^{9} = X^{9}$ and $V^{9}
= Y^{9}$. If $U^{6} = X^{6}$ and $V^{6} = Y^{6}$, then of course
we have $U^{3} = X^{3}$ and $V^{3} = Y^{3}$ and the lemma follows.
If $U^{6} = Y^{6}$ and $V^{6} = X^{6}$, then $U$, $V$, $X$ and $Y$
have the same norm. However, since $U + V = X + Y$, the pairs
$\{U, V\}$ and $\{X, Y\}$ are the same. hence implying the lemma.
The reason for this comes from the following statement which is
elementary: (Note that even when $U + V = X + Y = 0$, although we
cannot directly apply this statement, we still have $U^{3} + V^{3}
= X^{3} + Y^{3} = 0$, so that $\Set{U^{3}, V^{3}} = \Set{X^{3},
Y^{3}}$.)

\medskip

{\bf Statement} \emph{The pair $(z_{1}, z_{2})$ such that $|z_{1}|
= |z_{2}| = R$ and $z_{1} + z_{2} = Z$ is uniquely determined by
$R > 0$ and $Z$ with $0 < |Z| < 2 R$. }
\bigskip

So in all cases, $\Set{U^{3}, V^{3}} = \Set{X^{3}, Y^{3}}$.

\qedsymbol

\bigskip

\emph{Proof of Lemma ~\ref{T:703}} (contd.)

By the previous lemma,
\[
    \Set{U_{v}^{3}, V_{v}^{3}}
    = \Set{X_{v}^{3}, Y_{v}^{3}}
\]
at any unramified finite place $v$.

Hence
\begin{align}
    L &(s, {\rm sym}^{3} (\pi'_{v}))
    \notag \\
    &= {(1 - U_{v}^{3} {(Nv)}^{-s})}^{-1}
    {(1 - V_{v}^{3} {(Nv)}^{-s})}^{-1}
    \notag \\
    & {(1 - V_{v}^{2} U_{v} {(Nv)}^{-s})}^{-1}
    {(1 - U_{v}^{2} V_{v} {(Nv)}^{-s})}^{-1}
    \notag \\
    &= {(1 - X_{v}^{3} {(Nv)}^{-s})}^{-1}
    {(1 - Y_{v}^{3} {(Nv)}^{-s})}^{-1}
    \notag \\
    &{(1 - V_{v} \omega_{\pi'_{v}} ({\rm Frob}_{v}) {(Nv)}^{-s})}^{-1}
    {(1 - U_{v} \omega_{\pi'_{v}} ({\rm Frob}_{v}) {(Nv)}^{-s})}^{-1}
    \notag \\
    &= L (s, \chi_{1, v}^{3}) L (s, \chi_{2, v}^{3})
    L (s, \pi' \otimes \omega_{\pi'_{v}})
    \notag
\end{align}
Here we have used
\[
    U_{v} V_{v} = \omega_{\pi'_{v}} ({\rm Frob}_{v})
\]

Hence
\[
    {\rm sym}^{3} (\pi'_{v}) \, \simeq \, \chi_{1, v}^{3} \boxplus \chi_{2, v}^{3}
    \boxplus
    (\pi'_{v} \otimes \omega_{\pi'_{v}}),
\]

and by Lemma ~\ref{T:602} we get what we desire, namely,
\[
    {\rm sym}^{3} (\pi') \, \simeq \, \chi_{1}^{3} \boxplus \chi_{2}^{3}
    \boxplus (\pi' \otimes \omega_{\pi'}).
\]

\qedsymbol

\bigskip

\begin{lemma} \label{T:704}
    If $\pi'$ is tetrahedral with ${\rm sym}^{2} (\pi')$
    invariant under twisting by a cubic character $\chi$,
    then
    $$
    {\rm sym}^{3} (\pi') \, \cong \, (\pi' \otimes \omega_{\pi'}
    \chi)
    \boxplus (\pi' \otimes \omega_{\pi'} \chi^{-1}).
    $$
\end{lemma}

\medskip

Hence the situation of Lemma ~\ref{T:703} will not happen if
$\pi'$ is tetrahedral.

\medskip

\emph{Proof of Lemma ~\ref{T:704}.}

Consider $\pi' \boxtimes {\rm sym}^{2} (\pi') = {\rm sym}^{3}
(\pi') \boxplus (\pi' \otimes \omega_{\pi'})$, which obviously
contains $\pi' \otimes \omega_{\pi'}$ as an isobaric constituent .

Since ${\rm sym}^{2} (\pi')$ allows self twists by $\chi$ and
$\chi^{-1}$, the isobaric sum above should also contain $\pi'
\otimes \omega_{\pi'} \chi$ and $\pi' \otimes \omega_{\pi'}
\chi^{-1}$. Together with $\pi' \otimes \omega_{\pi'}$, they are
pairwise inequivalent, the reason being that if a cusp form on $GL
(2)$ admits a self twist by a character, then such character has
to be trivial or quadratic.

Thus, by the uniqueness of the isobaric decomposition, ${\rm
sym}^{3} (\pi')$ should have $\pi' \otimes \omega_{\pi'} \chi$ and
$\pi' \otimes \omega_{\pi'} \chi^{-1}$ as its constituents, and
there is no room for any other constituent.

\qedsymbol

\bigskip

\emph{Proof of Case III when $\pi'$ is not dihedral.}

\medskip

Assume that $\Pi = \pi \boxtimes \pi' = \sigma_{1} \boxplus
\sigma_{2}$ where $\sigma_{1}$ and $\sigma_{2}$ are cusp forms on
$GL (3)$ with central characters $\eta_{1}$ and $\eta_{2}$
respectively. Also, assume that $\pi'$ is not dihedral.

\medskip

\emph{Subcase A: \, $\pi$  does not allow a self twist by a
nontrivial character.}

\medskip

{F}rom \eqref{EQ:701} (Proposition ~\ref{T:701}) and
\eqref{EQ:703} (Lemma ~\ref{T:702}), we have

\begin{align}
    L (s, {\rm sym}^{3} &(\pi') \otimes \chi)
    L (s, (\pi \boxtimes \pi') \times \tilde{\pi} \otimes \omega_{\pi'} \chi)
    \notag \\
    &= L (s, \pi' \times \pi; \Lambda^{3} \otimes \omega_{\pi}^{-1} \chi)
    L (s, \pi' \otimes \omega_{\pi'} \chi)
    \notag \\
    &= L (s, \eta_{1} \omega_{\pi}^{-1} \chi) L (s, \eta_{2} \omega_{\pi}^{-1} \chi)
    L (s, \pi' \otimes \omega_{\pi'} \chi)
    \notag \\ \label{EQ:705}
    &L (s, \sigma_{1} \times \tilde{\sigma}_{2} \otimes \eta_{2} \omega_{\pi}^{-1} \chi)
    L (s, \sigma_{2} \times \tilde{\sigma}_{1} \otimes \eta_{1} \omega_{\pi}^{-1} \chi)
\end{align}

\medskip

Hence, take $\chi = \omega_{\pi} \eta_{i}^{-1}$, then the right
hand side has a pole at $s = 1$ (as the remaining factors do not
vanish st $s = 1$). Then the left hand side also has a pole at $s
= 1$.

\medskip

However, since $\pi'$ is not dihedral, then ${\rm sym}^{3} (\pi')$
is either cuspidal or an isobaric sum of two cusp forms on $GL
(2)$ (Lemma ~\ref{T:704}), so any twisted $L$--function of ${\rm
sym}^{3} (\pi')$ has to be entire. So the only pole at $s = 1$
should come from $L (s, (\pi \boxtimes \pi') \times \tilde{\pi}
\otimes \omega_{\pi'} \omega_{\pi} \eta_{i}^{-1})$.

\medskip

As $\pi$ is cuspidal, $\Pi = \pi \boxtimes \pi'$ should have both
$\sigma'_{i} = \pi \otimes \omega_{\pi'}^{-1} \omega_{\pi}^{-1}
\eta_{i}$ as constituents.

\medskip

Since $\pi$ is not monomial, it does not allow a self twist.
Hence, if $\eta_{1} \ne \eta_{2}$ then $\sigma'_{i}$ are
different. Hence $\sigma'_{1}$ and $\sigma'_{2}$ are the only
constituents of $\Pi$ which are also twists of $\pi$. Furthermore,
if $\eta_{1} = \eta_{2}$, then the order of the pole of both sides
of \eqref{EQ:705}, and hence also of $L (s, (\pi \boxtimes \pi')
\times \tilde{\pi} \otimes \omega_{\pi'} \omega_{\pi}
\eta_{i}^{-1})$, is $2$. Hence $\sigma'_{1} = \sigma_{2}$ should
be a constituent of $\Pi$ with multiplicity $2$.

\medskip

Thus we get an isobaric decomposition of $\Pi$ as a sum of two
twists of $\pi$. Thus, from lemma ~\ref{T:703} and Lemma
~\ref{T:704}, this cannot happen if $\pi'$ is not dihedral. This
completes the treatment of Subcase A.

\bigskip

\emph{Subcase B: $\pi$ admits a self twist by a nontrivial cubic
character $\chi$.}

\medskip

In this subcase, recall that we are assuming $\Pi = \pi \boxtimes
\pi' \cong \sigma_{1} \boxplus \sigma_{2}$, where $\sigma_{1}$ and
$\sigma_{2}$ are cusp forms with respective central characters
$\eta_{1}$ and $\eta_{2}$.

\medskip

We claim that $\sigma_{1}$ and $\sigma_{2}$ are also invariant
when twisted by $\chi$. Otherwise $\sigma_{i} \otimes \chi$ and
$\sigma_{i} \otimes \chi^{-1}$ will be distinct from $\sigma_{i}$,
while they should both be constituents of $\Pi \simeq \Pi \otimes
\chi \simeq \Pi$. Hence $\Pi$ has at least degree $3 \times 3 =
9$, which is impossible as it is automorphic on $GL (6)$.

\medskip

Moreover, let $v$ be any place where $\pi$ and $\pi'$ are
unramified. Write (for $i=1,2$) $\pi'_{v} = \alpha_{v} \boxplus
\alpha'_{v}$, $\pi_{v} = \beta_{v} \boxplus \beta_{v} \chi_{v}
\boxplus \beta_{v} \chi_{v}^{-1}$ (this form being implied by $\pi
\cong \pi \otimes \chi$), and $\sigma_{i, v} = \theta_{i, v}
\boxplus \theta_{i, v} \chi_{v} \boxplus \theta_{i, v}
\chi_{v}^{-1}$.

Then we have
\[
    \Pi_{v} = (\alpha_{v} \boxplus \alpha'_{v}) \otimes
    \beta_{v}\otimes
 (1 \boxplus \chi_{v} \boxplus \chi_{v}^{-1})
\]
and
\[
    \sigma_{1, v} \boxplus \sigma_{2, v}
    = (\theta_{1, v} \boxplus \theta_{2, v})\otimes
 (1 \boxplus \chi_{v} \boxplus \chi_{v}^{-1})
\]

Since the sets of all cubes of characters occurring in the
previous two isobaric decompositions should be the same, and since
$\beta_{v}^{3} = \omega_{\pi_{v}}$ and $\theta_{v}^{3} = \eta_{i,
v}$, we must have
\begin{align}
    \Set{\alpha_{v}^{3} \omega_{\pi_{v}}, {\alpha'_{v}}^{3} \omega_{\pi_{v}}}
    &= \Set{\alpha_{v}^{3} \beta_{v}^{3}, {\alpha'_{v}}^{3} \beta_{v}^{3}}
    \notag \\
    &= \Set{\theta_{1, v}^{3}, \theta_{2, v}^{3}}
    \notag \\
    &= \Set{\eta_{1, v}, \eta_{2, v}}
    \notag
\end{align}

So
\begin{align}
    {\rm sym}^{3} &(\pi'_{v}) \otimes \omega_{\pi_{v}}
    \cong \alpha_{v}^{3} \omega_{\pi_{v}} \boxplus {\alpha'_{v}}^{3} \omega_{\pi_{v}}
    \boxplus (\alpha_{v} \boxplus \alpha'_{v}) \boxtimes\omega_{\pi'_{v}} \omega_{\pi_{v}}
    \notag \\
    &\simeq \eta_{1, v} \boxplus \eta_{2, v}
    \boxplus \pi'_{v} \otimes \omega_{\pi'_{v}} \omega_{\pi_{v}},
    \notag
\end{align}
that is
\[
    {\rm sym}^{3} (\pi'_{v}) \cong
    \eta_{1, v} \omega_{\pi_{v}}^{-1} \boxplus
    \eta_{2, v} \omega_{\pi_{v}}^{-1}
    \boxplus \pi'_{v} \otimes \omega_{\pi'_{v}}
\]

\medskip

Thus by Lemma ~\ref{T:602} or the strong multiplicity one theorem,
we have
\[
    {\rm sym}^{3} (\pi')
    \cong \eta_{1} \omega_{\pi}^{-1} \boxplus \eta_{2} \omega_{\pi}^{-1}
    \boxplus  \pi' \otimes \omega_{\pi'}
\]

However, since $\pi'$ is not dihedral, ${\rm sym}^{3} (\pi')$ is
by Lemma ~\ref{T:704} either cuspidal or an isobaric sum of two
cusp forms on $GL (2)$. This gives a contradiction.

\medskip

This completes the treatment of Subcase B.

\medskip

\qedsymbol

\bigskip


\section{Proof of Theorem A, Part \#3}

In this part, we will treat the case when $\pi'$ is dihedral.
After that, we will analyze precisely the cuspidality criterion
when $\pi$ is an adjoint of a form on GL$(2)$. Again, $F$ denotes
a number field.

\medskip

In fact, we will prove the following:

\begin{theorem} \label{T:801}
    Let $\pi', \pi$ be cusp forms on $GL (2)/F, GL (3)/F$
    respectively, with $\pi'$ dihedral.

    Then $\Pi = \pi \boxtimes \pi'$ is cuspidal unless
    both the following two conditions hold:

    \textnormal{(a)} There is a non-normal cubic
    extension $K'$ of $F$ such that $\pi'_{K'}$ is
    Eisensteinian; equivalently, $\pi'$ is dihedral of type $D_{6}$.

    \textnormal{(b)} $\pi$ is monomial and $\pi = I^{F}_{K'} (\chi')$
    for some idele class character $\chi'$ of $K$.

    If \textnormal{(a)} and \textnormal{(b)} both hold,
    then $\Pi$ is an isobaric sum of two cuspidal
    representations of degree $3$, which are both
    twist equivalent to $\pi$.
\end{theorem}

\bigskip

Before we prove this theorem, let us recall that a dihedral Galois
representation $\rho'$ of $\G_{F}$ is said to be of type $D_{2 n}$
if its projective image is $D_{2 n}$. It is clear that $\rho'$ is
not irreducible if and only if $n = 1$ (Note that the projective
image of $D_{4 n}$ must be a quotient of $D_{2 n}$ since $D_{4 n}$
has a nontrivial center). If $6 | n$, then $D_{2 n}$ has a unique
cyclic subgroup with quotient isomorphic to $D_{6} \cong S_{3}$.
Suppose $K'$ is a non-normal cubic extension of $F$, and $\rho'$
restricted to $\G_{K'}$ is reducible. Then the projective image of
$\G_{K'}$ should be a subgroup of that of $\G_{F}$ of index $3$,
hence is isomorphic to $D_{2 n / 3}$. Thus $n = 3$, and $\rho$
must be dihedral of type $D_{6}$. Similarly, we conclude that if
$\pi'$ is dihedral and $\pi_{K'}$ is not cuspidal, then $\pi'$ is
of type $D_{6}$.

\bigskip

\emph{Proof of Theorem ~\ref{T:801}.}

First assume (a) and (b). Note that
$\pi'_{K'} = v_{1} \boxplus v_{2}$,
and $\pi = I^{F}_{K'} (\chi')$.
Then
\begin{align}
    \Pi &= \pi \boxtimes \pi' \cong \pi' \boxtimes I^{F}_{K'} (\chi')
    \notag \\
    &\cong I^{F}_{K'} (\pi'_{K'} \otimes \chi')
    \cong I^{F}_{K'} (v_{1} \chi') \boxplus I^{F}_{K'} (v_{2} \chi')
    \notag
\end{align}
Hence $\Pi$ is not cuspidal.

\medskip

Note that $\pi_{i} = I^{F}_{K'} (v_{i} \chi')$ MUST be cuspidal
as from Section 5, $\Pi$ cannot have a character as its constituent.

\medskip

Next, we prove that if $\Pi$ is not cuspidal,
then (a), (b) and the remaining statements
of the theorem hold.

\medskip

\emph{Step 1: $\pi_{K}$ is cyclic cubic monomial.}

\medskip

Assume that $\pi' = I^{F}_{K} (\tau)$ where $\tau$ is some idele
class character of $C_{K}$ with $K$ a quadratic extension of $F$.
and also assume that $K / F$ is cut out by $\delta$.

\medskip

{F}rom Section 5, Case I and II cannot happen, so we are in Case
III. Say $\Pi = \sigma_{1} \boxplus \sigma_{2}$ where $\sigma_{i}$
are some cusp forms on $GL (3) / F$. As $\pi'$ allows a self twist
by $\delta$, so does $\Pi = \pi \boxtimes \pi'$. Thus $\sigma_{1}
\cong \sigma_{2} \otimes \delta$ as the only possible characters
that either $\sigma_{i}$ allows (for self-twisting) should be
trivial or cubic.

\medskip

Let $\theta$ be the nontrivial Galois conjugation of $K / F$. Then
$\pi'_{K} \cong \tau \boxplus \tau^{\theta}$. Hence the base
change $\Pi_{K}= \pi'_{K} \boxtimes \pi_{K}$ is equivalent to
$\pi_{K} \otimes \tau$ plus $\pi_{K} \otimes (\tau\circ{\theta})$.
As $\Pi = \sigma_{1} \boxplus \sigma_{1} \otimes \delta$ hence
$\Pi_{K}$ is equivalent to the isobaric sum of two copies of
${\sigma_{1}}_K$.

Thus $\pi_{K} \otimes \tau \cong \pi_{K} \otimes
(\tau\circ{\theta}) \cong \sigma_{1}$. As $\pi' = I^{F}_{K}
(\tau)$ is cuspidal, $\tau \ne \tau\circ{\theta}$. Hence $\pi_{K}$
is forced to be cyclic monomial.

\medskip

\emph{Step 2: $\pi$ is non-normal cubic monomial.}

\medskip

By Step 1, $\mu = \tau^{-1} (\tau\circ{\theta})$ is a cubic
character of $C_{K}$. Let $M$ be the cubic field extension of $K$
cut out by $\mu$. As $\mu\circ{\theta} = \mu^{-1}$, $M^{\theta} =
M$, thus $M / F$ is normal, and ${\rm Gal}(M / F) \cong S_{3}$.

\medskip

Furthermore, $\pi_{K} = I^{K}_{M} (\lambda)$ for some character
$\lambda$ of $C_{M}$. And also, $\pi_{M}$ is of the form $\lambda
\boxplus \lambda' \boxplus \lambda''$.

\medskip

Let $K'$ be a non-normal cubic extension of
$F$ contained in $M$. Then $[M : K'] = 2$
and $\pi_{M}$ is a quadratic base change of $\pi_{K'}$.

We claim that $\pi_{K'}$ is Eisensteinian, i.e., not cuspidal. The
reason is that, if $\pi_{K'}$ were cuspidal, then its quadratic
base change $\pi_{M}$ would be either cuspidal or the isobaric sum
of two cusp forms of the same degree.  Since $\pi_{K'}$ is a cusp
form on $GL (3)$, we see  from \cite{AC} that this is impossible.

\medskip

So $\pi_{K'}$ must admit a character as an isobaric constituent,
which means that $\pi$ is induced from some character of $C_{K'}$.

\medskip

\emph{Step 3: $\pi'_{K'}$ is not cuspidal, hence $\pi'$ is dihedral of type
$D_{6}$.}

\medskip

Recall that $\pi'_{K} = \tau \boxplus (\tau\circ{\theta}) = \tau
\boxplus \tau \mu$, so that $\pi'_{M} = \tau_{M} \boxplus
\tau_{M}$ as $M / K$ is cut out by $\mu$.

\medskip

Thus the projective image of $\rho'_{M}$ is trivial, where $\rho'$
is the representation Ind$_K^F(\tau)$ of the Weil group $W_F$, and
$\rho'_{M}$ is the restriction of $\rho$ to $\G_{M}$. Hence the
projective image of $\rho$ must be $D_{6}$.

\medskip

\emph{Remark: } Even if $\pi'$ is selfdual, $\tau$ may be a
character of order $3$ or $6$.

\medskip

\emph{Step 4: $\sigma_{1}$ and $\sigma_{2}$ are all twist equivalent to $\pi$.}

\medskip

Observe that
\begin{align}
    {(\tau \mu^{-1})}\circ{\theta} &= ({\tau}\circ{\theta}) ({\mu}\circ{\theta})^{-1}
    \notag \\
    &= \tau \mu \mu = \tau \mu^{-1}.
    \notag
\end{align}

So $\tau \mu^{-1}$ is a base change of some character, say $\nu$,
of $C_{F}$ to $K / F$.

\medskip

So $\pi' = I^{F}_{K} (\mu) \otimes \nu$ and $\pi'_{K'} \cong
\nu_{K'} \boxplus \nu_{K'} \delta_{K'}$. By Step 2, we may assume
that $\pi = I^{F}_{K'} (\lambda)$ for a character $\lambda$ of a
non-normal cubic extension $K'$ of $F$. We get
\begin{align}
    \Pi &= \pi \boxtimes \pi' = \pi' \boxtimes I^{F}_{K'} (\lambda)
    \notag \\
    &= I^{F}_{K'} (\pi'_{K'} \otimes \lambda)
    = I^{F}_{K'} (\nu_{K'} \lambda \boxplus \nu_{K'} \lambda \delta_{K'})
    \notag \\
    &= I^{F}_{K'} (\lambda) \otimes \nu \boxplus I^{F}_{K'} (\lambda)
    \otimes \nu \delta
    \notag \\
    &= (\pi \otimes \nu) \boxplus (\pi \otimes \nu \delta)
\end{align}

\medskip

Now the proof of Theorem 8.1 is completed.

\qedsymbol

\bigskip

\emph{Remark. } When $\pi$ is twist equivalent to $Ad (\pi_{0})$,
and $\pi'$ is dihedral, we claim that the only way $\Pi = \pi
\boxtimes \pi'$ can be {\it not} cuspidal is for $Ad (\pi)$ to be
non-normal cubic monomial, implying that $\pi$ is of octahedral
type. We get the following theorem which is more precise than the
result in \cite{Wa2003}:

\bigskip

\begin{theorem} \label{T:802}
    Let $\pi'$, and $\pi''$ be two cusp form on $GL (2) / F$,
    and suppose that $\pi''$ is not dihedral.
    then $\Pi = \pi' \boxtimes Ad (\pi'')$ is cuspidal
    unless one of the following happens:

    \textnormal{(d)} $\pi'$ and $\pi''$ are twist equivalent.

    \textnormal{(e)} $\pi'$ and $\pi''$ are octahedral attached with
    the same $\tilde{S}_{4}$-extension, and $Ad (\pi')$
    and $Ad (\pi'')$ are twist equivalent.

    \textnormal{(f)} $\pi''$ is octahedral, $\pi' = I^{F}_{K} (\mu) \otimes \nu$
    is dihedral, where $\mu$ is the cubic character which is allowed by
    $Ad (\pi''_{K})$.
\end{theorem}

\bigskip

The proof of the theorem above uses the following proposition.

\begin{proposition} \label{T:803}
    Let $\pi_{1}$ and $\pi_{2}$ are two non-dihedral cusp forms
    on $GL (2) / F$, and $Ad (\pi_{1})$ and $Ad (\pi_{2})$
    are twist equivalent. Then one of the following holds.

    \textnormal{(g)} $\pi_{1}$ and $\pi_{2}$ are twist equivalent
    (so that their adjoints are equivalent).

    \textnormal{(h)} $\pi_{1}$ and $\pi_{2}$ are octahedral attached with
    the same $\tilde{S}_{4}$-extension, and $Ad (\pi_{1})$
    and $Ad (\pi_{2})$ are twist equivalent by a quadratic characters.
\end{proposition}

\bigskip

\emph{Proof of Proposition ~\ref{T:803}.} (cf.\ \cite{Ra2000})

It is clear that (g) and (h) imply that $\pi_{1}$ and $\pi_{2}$
are twist equivalent. So it suffices to show the other side.

\medskip

First assume that $Ad (\pi_{1})$ and $Ad (\pi_{2})$ are
equivalent.

Consider $\Pi = \pi_{1} \boxtimes \pi_{2}$.
Note that
$$
\Pi \boxtimes \overline{\Pi} \cong 1 \boxplus Ad (\pi_{1})
\boxplus Ad (\pi_{2}) \boxplus Ad (\pi_{1}) \boxtimes Ad (\pi_{2})
$$
admits two copies $1$ as its constituents. Hence $\Pi$ is not
cuspidal.

\medskip

If $\Pi$ contains a character $\nu$, then $\overline{\pi}_{2}
\cong \pi_{1} \otimes \nu$. If $\Pi$ is an isobaric sum of two
cusp forms $\sigma_{1}$ and $\sigma_{2}$ on $GL (2)$, then check
that $\Lambda^{2} (\pi_{1} \boxtimes \pi_{2})$ is equivalent to
$$
(Ad (\pi_{1}) \boxplus
Ad (\pi_{2})) \otimes \omega_{\pi_{1}} \omega_{\pi_{2}}
$$
which does not contain any character constituent;
However, $\Lambda^{2} (\sigma_{1} \boxplus \sigma_{2})$
is equivalent
$$
\omega_{\sigma_{1}} \boxplus \omega_{\sigma_{1}} \boxplus
(\sigma_{1} \boxtimes \sigma_{2})
$$
which contains two $GL (1)$-constituents. Thus we get a
contradiction, and (g) holds.

\medskip

Furthermore, assume that
$$Ad (\pi_{2}) \cong Ad (\pi_{1}) \otimes \epsilon$$
where $\epsilon$ is a character. Then
$$Ad (\pi_{2}) \cong Ad (\pi_{1}) \otimes \epsilon^{-1}$$
and hence
\[
    Ad (\pi_{1}) \boxtimes Ad (\pi_{1})
    \cong Ad (\pi_{2}) \boxtimes Ad (\pi_{2}).
\]

\medskip

However,
\[
    Ad (\pi_{i}) \boxtimes Ad (\pi_{i}) \cong 1
    \boxplus Ad (\pi_{i}) \boxplus A^{4} (\pi_{i}).
\]
Hence
\[
    Ad (\pi_{1}) \boxplus A^{4} (\pi_{1})
    \cong Ad (\pi_{2}) \boxplus A^{4} (\pi_{2}).
\]

If $Ad (\pi_{1})$ and $Ad (\pi_{2})$ are equivalent,
we get (g). Otherwise, $Ad (\pi_{1})$,
which is a nontrivial twist of $Ad (\pi_{2})$,
must be contained in $A^{4} (\pi_{2})$.

\medskip

So $\pi_{1}$ and $\pi_{2}$ are of solvable polyhedral type.

If $\pi_{2}$ is tetrahedral, then
\[
    A^{4} (\pi_{2}) \cong Ad (\pi_{2}) \boxplus \omega \boxplus \omega^{2}
\]
where $\omega$ is some cubic character that $Ad (\pi_{2})$ admits
as as self-twist. So this cannot happen.

\medskip

Thus $\pi_{2}$ and $\pi_{1}$
are octahedral, and
\[
    A^{4} (\pi_{2}) \cong I^{F}_{K} (\mu) \boxplus Ad (\pi_{2})
    \otimes \epsilon
\]
where $K$ is a quadratic extension
of $F$ such that $Ad ({\pi_{2}}_{K})$
allows a self twist by $\mu$,
and $\epsilon$ is the quadratic character cuts out $K$.

\medskip

So they must come from the same $\tilde{S}_{4}$-extension, and
hence (h) holds.

\qedsymbol

\bigskip

\emph{Proof of Theorem ~\ref{T:802}.}

Set $\pi = Ad (\pi'')$ {F}rom what we have seen (including the
proof of Theorem ~\ref{T:801}), $\Pi = \pi' \boxtimes Ad (\pi')$
is cuspidal unless (i) $Ad (\pi')$ and $Ad (\pi'')$ are twist
equivalent; or (ii) $\pi'$ is dihedral of type $D_{6}$, $\pi = Ad
(\pi'')$ is non-normal cubic monomial, and (a) and (b) of Theorem
~\ref{T:801} hold.

\medskip

If (ii) holds, then $\pi''$ must be octahedral, and (f) must hold.
If (1) holds, then $Ad (\pi')$ and $Ad (\pi'')$ are twist
equivalent. Then part (g) of Proposition ~\ref{T:803} implies (d),
and part (h) of this proposition implies (e).

\qedsymbol

\medskip

The proof of Theorem A is now complete.

\bigskip


\section{Proof of Theorem B}

In this part we deduce Theorem ~\ref{TM:B} from Theorem
~\ref{TM:A}. First we need some preliminaries.

\medskip

Recall that a cusp form $\pi$ on $GL (n)$ over $F$ is essentially
selfdual if $\overline{\pi}$ is twist equivalent to $\pi$.
Throughout this section, $\pi'$ and $\pi$ denote cusp forms on $GL
(2)$ and $GL (3)$ over $F$. We assume that $\pi'$ is not dihedral,
and $\pi$ is not twist equivalent to $Ad (\pi'')$ for any cusp
form $\pi''$.

\medskip

First, from Theorem ~\ref{TM:A},
we have the following:

\begin{corollary} \label{T:901}
    If $\pi'$ is not of solvable polyhedral type
    and $\pi$ is not essentially selfdual,
    then $\Pi = \pi \boxtimes \pi'$ is cuspidal.
\end{corollary}

\bigskip

\begin{lemma} \label{T:902}
    Let $K$ be any solvable extension of $F$.
    If $\pi$ is not essentially selfdual, and
    if $\pi_{K}$ does not admit any self twist,
    then $\pi_{K}$ is not essentially selfdual.
\end{lemma}

\medskip

\emph{Proof of Lemma ~\ref{T:902}.}

First, assume that $[K : F] = l$ a prime
so that $K / F$ is cyclic.
Let $\theta \ne 1$ be a Galois conjugation of
$K$ over $F$, and $\tau$ be a character
cutting out $K / F$.

\medskip

Assume that $\pi_{K}$ is essentially selfdual. Say
$\overline{\pi}_{K} \cong \pi_{K} \otimes \mu$, for a character
$\mu$. Applying $\theta$, and being aware of the fact that
$\pi_{K}$ is fixed by $\theta$, we get
\[
    \overline{\pi}_{K} \cong \pi_{K} \otimes (\mu\circ{\theta})
\]

Since $\pi_{K}$ does not allow a self twist, then
$\mu\circ{\theta} = \mu$, hence $\mu$ must be a base change of
some character $\alpha$ of $C_{F}$ to $K$.

\medskip

Hence, $\overline{\pi}$ and $\pi \otimes \alpha$ have the same
base change over $K / F$, and thus must be twist equivalent. This
shows that $\pi$ is also essentially selfdual.

\medskip

In general case, let $K_{0} = F
 \subset K_{1} \subset \ldots \subset K_{n} = K$
be a tower of cyclic extensions of prime degree.
Assume that $\pi_{K_{n}} = \pi_{K}$
is essentially selfdual,
then as $\pi_{K_{n}}$ does not allow a self twist,
neither does $\pi_{K_{i}}$ for any smaller $i$,
thus applying the arguments above inductively,
we claim that $\pi_{K_{i}}$
is essentially selfdual. In particular,
$\pi$ must be essentially selfdual.

\qedsymbol

\bigskip

\emph{Proof of Theorem ~\ref{TM:B}}.

First prove (a). $\Pi$ is cuspidal from Corollary ~\ref{T:901}.
First, assume only that $\pi$ is not essentially selfdual and does
not allow a self twist. Assume that $\Pi$ allows a
self twist by some character $\nu$.
Without loss of generality, may assume that
$\nu$ is of order $2$ or $3$. Let $K / F$
be cut out by $\nu$. Thus $\Pi_{K} = I^{F}_{K} (\eta)$
is Eisensteinian of type $(2, 2, 2)$
or $(3, 3)$.

However, $\pi_{K}$ is cuspidal as $\pi$ does not allow a self
twist. From Theorem ~\ref{TM:A} (and the remark at the end of
Section 5), $\Pi_{K}$ cannot be of type $(3, 3)$ as $\pi'_{K}$ is
not dihedral, type $(2, 2, 2)$ as $\pi'_{K}$ is not tetrahedral
(as $\pi'$ and hence $\pi'_{K}$ is not of solvable polyhedral
type). Thus $\Pi_{K}$ must be cuspidal, and hence $\Pi$ does not
allow a self twist.

\medskip

Moreover, assume that $\pi$ is not monomial,
in particular, $\pi$ is not induced from
a non-normal cubic extension.
Want to prove that $\Pi$ is not either.

\medskip

Assume that $\Pi \cong I^{F}_{K'} (\eta)$ where $\eta$ is some
cusp form on $GL (2)$ over $K'$ which is non-normal cubic over
$F$. Let $M$ be the normal closure of $K'$ over $F$ and $E$ be the
unique quadratic subextension of $M$ over $F$. Then $\pi'_{E}$ is
still not of solvable polyhedral type. And $\pi_{E}$ is not cyclic
monomial as $\pi$ is not monomial. From Lemma ~\ref{T:902},
$\pi_{E}$ is not essentially selfdual. Thus the first part of (a)
implies that $\Pi_{E}$ does not allow a self twist. Hence
$\Pi_{M}$ is still cuspidal. However, $\Pi_{M} \cong {I^{F}_{K'}
(\eta)}_{M} = \eta_{M} \boxplus (\eta_{M}\circ{\theta'}) \boxplus
(\eta_{M}\circ{\theta'^{2}})$ where $\theta'$ is the character
cutting out $M / E$. We get a contradiction.

Thus (a) is proved.

\bigskip

Next prove (b). It suffices to prove that
$\Pi_{K}$ is cuspidal for any solvable
extension $K$ of $F$.

Since $\pi'$ is not of solvable polyhedral type, neither is
$\pi'_{K}$. As $\pi$ is not of solvable type, then $\pi_{K}$ must
be cuspidal. We claim that $\pi_{K}$ cannot allow a self twist.
Otherwise, say $\pi_{K} \cong \pi_{K} \otimes \nu$. $\nu$ must be
a cubic character. Let $K_{1} / K$ be cut out by $\nu$, then
$\pi_{K_{1}}$ should be Eisensteinian. However, $K_{1} / F$ is
contained in some solvable normal extension. Thus $\pi$ is of
solvable type. Contradiction.

\medskip

Hence the claim. From Lemma ~\ref{T:902},
$\pi_{K}$ is not essentially selfdual.
Corollary ~\ref{T:901} then implies that
$\Pi_{K}$ is cuspidal. Thus $\Pi$
is not of solvable type.

\qedsymbol

\bigskip

Before we finish this section, we would like to
prove the following lemma.

\begin{lemma} \label{T:903}
    Let $\pi$ be a cusp form on $GL (2 m + 1)$ over $F$.
    Assume that $\pi$ is regular algebraic at infinity, and
    $F$ is not totally complex. Then $\pi$
    is not monomial.
\end{lemma}

\emph{Proof of Lemma ~\ref{T:903}.}

Assume that $\pi = I^{F}_{K'} (\nu)$ where $K'$ is a field
extension of $F$ of degree $2 m + 1$, and $\nu$ is an algebraic
character of $C_{K}$. As $F$ is not totally complex, neither is
$K'$ as $[K' : F]$ is odd. Thus from Weil (\cite{Weil}), $\nu$
must be of the form $\nu_{0} {\| \cdot \|}^{k}$ where $\nu_{0}$ is
a character of finite order. Thus $I^{F}_{K'} (\nu)$ does not
contain any nontrivial algebraic character at infinity, and hence
cannot be regular at infinity. \qedsymbol

\bigskip


\section*{Appendix}
\addtocounter{section}{1}
\setcounter{theorem}{0}

\bigskip

The object here is to justify the statement made in the
Introduction that it is possible to construct, for $n > 2$ even
(resp. $n=4$), non-selfdual, {\it monomial} (resp. non-monomial,
but imprimitive), cuspidal cohomology classes for $\Gamma \subset
{\rm SL}(n, \Z)$.

\medskip

\begin{main} \label{TM:E}
\begin{enumerate}
\item[(a)] \, Fix any integer $m > 1$. Then there exists a cuspidal
automorphic representation $\pi$ of GL$(2m, \A_\Q)$, which is not
essentially self-dual, contributing to the cuspidal cohomology in
degree $m^2$, and admitting a self-twist relative to a character
of order $2m$. In fact, this can be done for any, not necessarily
constant, coefficient system.
\item[(b)] \, There exists a cusp form
$\pi$ on GL$(4)/\Q$ contributing to the cuspidal cohomology in
degree $4$, which is not essentially self-dual. It admits a
self-twist relative to a quadratic character, but not relative to
any quartic character.
\end{enumerate}
\end{main}

\medskip

In both cases it will be apparent from the proof below that there
are infinitely many such examples, and by the discussion in
section 3, they give rise, for arbitrary coefficient systems $V$,
to non-selfual, cuspidal cohomology classes in $H^\ast(\Gamma, V)$
for suitable congruence subgroups $\Gamma \subset {\rm SL}(2m,
\Z)$ ($m \geq 2$). Moreover there are naturally associated Galois
representations, which are monomial in case (a), and are
imprimitive but non-monomial in case (b).

\bigskip

{\it Proof}. \, (a) \, For any $m > 1$, fix a finite-dimensional
coefficient system $V$. Then a cuspidal automorphic representation
$\pi$ of GL$(2m, \A_F)$ contributes to the cuspidal cohomology
with coefficients in $V$ iff it is algebraic with infinity type
(in the {\it unitary} normalization): (cf. \cite{C})
$$
\{(z/\vert z\vert)^{k_1}, ({\overline z}/\vert z\vert)^{k_1},
(z/\vert z\vert)^{k_2}, ({\overline z}/\vert z\vert)^{k_2}, \dots,
(z/\vert z\vert)^{k_m}, ({\overline z}/\vert z\vert)^{k_m}\},
$$
where $(k_1, k_2, \dots, k_m)$ is an ordered $m$-tuple of integers
(determined by $V$) satisfying
$$
k_1 > k_2 > \dots > k_m.
$$
In particular, $\pi_\infty$ is regular. If $V \simeq \C$, then as
seen in Theorem 3.2, $k_j = 2(m-j)+1$. Now pick any cyclic,
totally real extension $F$ of $\Q$ of degree $m$, and a totally
imaginary quadratic extension $K$ of $F$ which is normal over
$\Q$. Let $v_1, \dots, v_m$ denote the archimedean places of $F$,
and for each $j$ let $w_j$ be a complex place of $K$ above $v_j$.
Choose an algebraic Hecke character $\chi$ of of $K$ such that
$$
\chi_{w_j}(z) \, = \,
\left(\frac{z}{|z|}\right)^{2(m-j)+1}|z\overline z|^{2m-1} \, \,
\forall \, \, j \leq m. \leqno(1)
$$
That such a character exists is a consequence of the discussion on
page 3 of \cite{Weil}. To elaborate a bit for the sake of the
uninitiated, the necessary and sufficient condition for the
existence of $\chi$ as above is that there be a positive integer
$M$ such that the following holds for all units $u \in \mathfrak
O_K^\ast$ with components $u_j$ at $w_j$:
$$
\prod_{j \leq m} \, \left(\frac{u_j}{\overline u_j}\right)^{Mk_j}
\, = \, 1.\leqno(2)
$$
But the index of the real units $\mathfrak O_F^\ast$ in $\mathfrak
O_K^\ast$ is finite by the Dirichlet unit theorem, and hence for a
suitable $M$, $u_j^M$ is real for all $j$. The desired identity
(2) follows.

\medskip

Next pick a finite order character $\nu$ of $K$ and set
$$
\Psi \, = \, \chi\nu. \leqno(3)
$$
Let $\tau$ be the non-trivial automorphism ({\it complex
conjugation}) of $K/F$, and $\delta$ the quadratic character of
$F$ attached to $K$. Then we may, and we will, choose $\nu$ in
such a way that
$$
\Psi(\Psi\circ\tau) \, \ne \, \mu\circ N_{K/\Q}
\leqno(4)
$$
for any character $\mu$ of $\Q$, which is possible -- and this is
{\it crucial} -- because $[K:\Q] \geq 4$ and so $F = \{x \in K \,
\vert \, x^\tau = x\} \, \ne \, \Q$. Put
$$
\pi : = \, I_K^\Q(\Psi),
$$
where $I$ denotes {\it automorphic induction} (\cite{AC}). Note
that $\pi$ makes sense because $K/\Q$ is solvable and normal,
allowing automorphic induction to be defined. By looking at the
infinity type (1) we see that $\pi$ is regular and algebraic.

By construction, $\pi_\infty$ contributes to cohomology, and $\pi$
is cuspidal because the infinity type of $\Psi$ precludes it from
being $\Psi \circ \sigma$ for any non-trivial automorphism
$\sigma$ of $K/\Q$. To elaborate, note first that $\eta:=
I_K^F(\Psi)$ is cuspidal and algebraic, corresponding to a Hilbert
modular newform on GL$(2)/F$ of the prescribed weight at infinity.
Since the automorphic induction is compatible with doing it in
stages, $\pi$ is just $I_F^\Q(\eta)$, and since $F/\Q$ is cyclic,
it suffices to check that for any automorphism $\tau$ of $F$,
$\eta$ and $\eta\circ \tau$ are not isomorphic, which is clear
from the properties of $\Psi$.

It remains to check that $\pi$ is not essentially self-dual, which
comes down to checking the same for the $2m$-dimensional
representation $\rho$ of $W_\Q$ induced by the character $\Psi$ of
$W_K$. For this we need, by Mackey, to check that $\Psi^{-1} \ne
(\mu_K)\Psi\circ \sigma$ for any automorphism $\sigma$ of $K$ and
any character $\mu$ of $\Q$. By our choice of the infinity type,
this is clear for any $\sigma \ne \tau$ (and any $\mu$). For
$\sigma = \tau$, this is the content of (4). So we are done with
the proof of part (a).

\bigskip

(b) \, Let $K/\Q$ be an imaginary quadratic field and $\beta$ a
non-dihedral cusp form of weight $2$ over $K$ with $\Q$
coefficients, such that no twist of $\beta$ is a base change from
$\Q$. Here {\it weight 2} signifies the fact that the Langlands
parameter of $\beta_\infty$ is given by
$$
\sigma(\beta_\infty) \, = \, \{z/|z|, \overline z/|z|\} \otimes
(z\overline z)
$$
Here are two explicit (known) examples with these properties:
First consider the (non-CM) elliptic curve
$$
E_1: \, y^2+xy \, = \, x^3+(3+\sqrt{-3})x^2/2+(1+\sqrt{-3})x/2
$$
over $K_1=\Q(\sqrt{-3})$. This was shown to be associated to a
cusp form $\beta_1$ of weight $2$ and trivial central character on
GL$(2)/K_1$ by R.~Taylor (\cite{Ta}) such that $a_P(E_1) =
a_P(\beta_1)$ for a set of primes $P$ of density $1$. (In fact,
recent results can be used to show that this equality holds
outside a finite set of primes $P$.)

For the second example, set $K_2 = \Q(i)$ and $Q$ the prime ideal
generated by $8+13i$. Then there is a corresponding cusp form
$\beta_2$ of weight $2$, conductor $Q$  and trivial central
character, as seen on page 344 of the book \cite{EGM} by Elstrodt,
Grunewald and Mennicke. Its conjugate by the non-trivial
automorphism $\theta$ of $K_2$ will have conductor $Q^\theta$ and
so no twist of $\beta_2$ can be a base change from $\Q$. There is
a corresponding elliptic curve
$$
E_2: \, y^2 +iy \, = \, x^3 +(1+i)x^2 +ix
$$
over $K_2$ of conductor $Q$, and one knows for many $P$ that
$a_P(E_2) = a_P(\beta_2)$.

Next choose an algebraic Hecke character $\chi$ of $K$ such that
$$
\chi_{\infty}(z) \, = \, ({z}/{|z|})^2|z\overline z|^{2}.
$$
For example, we can choose $\chi$ to be the square of a Hecke
character associated to a CM elliptic curve. Now consider, for $j
= 1, 2$, the automorphic induction
$$
\pi_j : = \, I_K^\Q(\beta_j \otimes \chi).
$$
The infinity types chosen imply that either $\beta_j \otimes \chi$
is not isomorphic to its transform by the non-trivial automorphism
of $K_j$. So $\pi_j$ is a cusp form on GL$(4)/\Q$. It is
cohomological, as easily seen by its archimedean parameter. That
$\pi_j$ is not essentially self-dual is an immediate consequence
of the infinity types of $\chi$ and $\beta$. Finally, $\pi_j$
admits a quadratic self-twist, namely by the character of $\Q$
associated to $K_j$, but it admits no quartic self-twist as
$\beta_j$ is non-dihedral. We are now done.

\qed

\bigskip

Dinakar Ramakrishnan

253-37 Caltech

Pasadena, CA 91125, USA.

dinakar@caltech.edu

\bigskip

Song Wang

Department of Mathematics

Yale University

New Haven, CT 06520, USA.

song.wang@yale.edu

\bigskip

\stop


\bigskip

\begin{thebibliography}{99}

\bibitem[AC1989]{AC} J.~ Arthur and L.~ Clozel,
\emph{Base Change, and the Advanced Theory of the Trace Formula},
Annals of Math.\ Studies \textbf{120}, Princeton, 1989.

\bibitem[AGG1984]{A-G-G} A.~Ash, D.~Grayson, P.~Green,
\emph{Computations of cuspidal cohomology of congruence subgroups
of ${\rm SL}(3, Z)$}, Journal of Number Theory {\bf 19} (1984),
no. 3, 412--436.

\bibitem[BaR1994]{Ba-R} L.~Barthel and D.~Ramakrishnan,
\emph{A non-vanishing
result for twists of $L$-functions of GL$(n)$}, Duke Math. Journal
{\bf 74}, No. 3 (1994), 681--700.

\bibitem[BoJ1979]{BoJ} A. ~Borel and H. ~Jacquet,
\emph{Automorphic forms and automorphic representations} (with a
supplement "On the notion of an automorphic representation" by R.
P. Langlands), Proc. Sympos. Pure Math., XXXIII, Automorphic
forms, representations and $L$-functions (Proc. Sympos. Pure
Math., Oregon State Univ., Corvallis, Ore., 1977), Part 1,
189--207, Amer. Math. Soc., Providence, R.I., 1979.

\bibitem[BoW1980]{BoW} A. ~Borel and N.R. ~Wallach,
\emph{Continuous cohomology, discrete subgroups,
 and representations of reductive groups},
Annals of Mathematics Studies {\bf 94},
Princeton University Press, Princeton, N.J.;
University of Tokyo Press, Tokyo, 1980.

\bibitem[C1988]{C} L. ~Clozel,
\emph{Motifs et formes automorphes: Applications
du principe de fonctorialit\'{e}},
Automorphic forms, Shimura varieties, and $L$-functions,
Vol. I (Ann Arbor, MI, 1988), 77--159.

\bibitem[EGM1998]{EGM} J.~Elstrodt, F.~Grunewald and J.~Mennicke,
{\it Groups acting on hyperbolic space}, Springer Monographs in
Math., Berlin (1998).

\bibitem[vGKTV1997]{vG-K-T-V} B.~van Geemen, W.~van der Kallen, J.~Top
and A.~Verberkmoes, Hecke eigenforms in the cohomology of
Congruence subgroups of SL$(3, \Z)$, Experimental Mathematics {\bf
6}:2, 163--174 (1997).

\bibitem[vGT1994]{vG-T1} B.~van Geemen and J.~Top, A non-selfdual
automorphic representation of GL$_3$ and a Galois representation,
Inventiones Math. {\bf 117} (1994), no.3, 391--401.

\bibitem[vGT2000]{vG-T2} B.~van Geemen and J.~Top,  Modular forms
for ${\rm GL}(3)$ and Galois representations, in {\it Algorithmic
number theory} (Leiden, 2000), 333--346, Lecture Notes in Comput.
Sci. {\bf 1838}, Springer, Berlin (2000).

\bibitem[HaT2000]{HaT2000} M.~ Harris and R.~ Taylor,
\emph{On the geometry and cohomology of some simple Shimura varieties},
With an appendix by Vladimir G. Berkovich.
Annals of Mathematics Studies, \textbf{151}. Princeton University Press, Princeton, NJ, 2001.

\bibitem[He2000]{He2000} G.~ Henniart,
\emph{Une preuve simple des conjectures de Langlands pour GL (n) sur un corps $p$--adique},
Inventiones Math. \textbf{139} no. 2 (2000), 439--455 (French).

\bibitem[JPSS1983]{JPSS} H.~Jacquet, I.~Piatetski-Shapiro and
J.~Shalika, Rankin-Selberg convolutions, Amer. J. Math.{\bf 105},
367--464 (1983).

\bibitem[JS1981]{JS81} H.~ Jacquet and J.A.~ Shalika,
\emph{Euler products and classification of automorphic forms} I
and II, Amer.\ J.\ of Math. \textbf{103} (1981), 499--558 \&
777--815.

\bibitem[JS1990]{JS90} H.~ Jacquet and J.A.~ Shalika,
\emph{Rankin--Selberg convolutions: Archimedean theory} in
\emph{Piatetski--Shapiro Festschrift}, Isreal Math.\ conf.\ Proc.,
Part II, 125--207, The Weizmann Science Press of Isreal (1990).

\bibitem[K2001]{K2001} H.~ Kim,
\emph{Functoriality of the exterior square of $GL_{4}$ and the
symmetric fourth power of $GL_{2}$},  J. Amer. Math. Soc. {\bf 16}
(2003), no.1, 139--183.

\bibitem[KSh2002-1]{KSh2000} H.~ Kim and F.~ Shahidi,
\emph{Functorial products for $GL(2) \times GL(3)$ and the
symmetric cube for $GL(2)$} (with an appendix by C.J.~ Bushnell
and G.~ Henniart), Annals of Math (2) \textbf{155}, (2002), no.3,
837--893.

\bibitem[KSh2002-2]{KSh2001} H.~ Kim and F.~ Shahidi,
\emph{Cuspidality of symmetric powers with applications}, Duke
Journal of Math (\textbf{112}), no.1, 177--197.

\bibitem[Ku1980]{Ku} S. ~Kumaresan,
\emph{On the canonical $k$-types in the irreducible unitary
$g$-modules with nonzero relative cohomology}, Invent. Math. {\bf
59} (1980), no. 1, 1--11.

\bibitem[La1979-1]{La79-1} R.P.~ Langlands,
\emph{On the notion of an Automorphic Representation. A Supplement
to the preceeding paper}, Proceedings of Symposia in Pure
Mathematics (Corvallis), Vol \textbf{33} (1979), part 1, 203--207.

\bibitem[La1979-2]{La79-2} R.P.~ Langlands,
\emph{Automorphic Representations, Shimura Varieties, and Motives.
Ein M\"{a}rchen}, Proceedings of Symposia in Pure Mathematics
(Corvallis), Vol \textbf{33} (1979), part 2, 205--246.

\bibitem[La1980]{La80} R.P.~ Langlands,
\emph{Base change for $GL (2)$}, Annals of Math.\ Studies
\textbf{96}, Princeton, 1980.

\bibitem[MW1995]{MW} C.~Moeglin and J.-L.~Waldspurger, Spectral
decomposition and Eisenstein series: Une paraphrase de
l'\'Ecriture, Cambridge University Press (1995).

\bibitem[Ra2000]{Ra2000} D.~ Ramakrishnan,
\emph{Modularity of the Rankin--Selberg $L$--series, and the
multiplicity one for $SL (2)$}, Ann. Math. \textbf{150} (2000),
45--111.

\bibitem[RaWa2001]{RaWa2001} D. ~Ramakrishnan and S.~ Wang,
\emph{On the exceptional zeros of Rankin--Selberg $L$-Functions},
Compositio Mathematica \textbf{135} (2003), 211--244.

\bibitem[Sh1988]{Sh} F.~Shahidi, On the Ramanujan conjecture and
finiteness of poles for certain $L$-functions.  Ann. of Math. (2)
{\bf 127} (1988), no. 3, 547--584.

\bibitem[WSt2003]{WSt} W.~Stein, {\it Modular forms database},
see \texttt{http://modular.fas.harvard.edu/Tables/index.html}.

\bibitem[Ta1994]{Ta} R.~Taylor, $l$-adic representations associated to
modular forms over imaginary quadratic fields II,  Inventiones
Math. {\bf 116} (1994),  no. 1-3, 619--643.

\bibitem[VZ1984]{VZ} D.A. ~Vogan, G.J. ~Zuckerman,
\emph{Unitary representations with nonzero cohomology},
Compositio Math. {\bf 53} (1984), no. 1, 51--90.

\bibitem[Wa2003]{Wa2003} Song Wang,
\emph{On the Symmetric Powers of Cusp Forms on $GL (2)$ of
Icosahedral Type}, IMRN {\bf 2003}:44 (2003), 2373--2390.

\bibitem[We1955]{Weil} A.~Weil,
\emph{On a certain type of characters of the idèle-class group of
an algebraic number-field}, Proceedings of the International
Symposium on Algebraic Number Theory, Tokyo \& Nikko, 1955, pp.
1--7, Science Council of Japan, Tokyo, 1956.

\end{thebibliography}
\end{document}